\documentclass[3p,final]{elsarticle}

\usepackage{amssymb}
\usepackage{amsmath}
\usepackage{braket}
\usepackage{mathtools}
\usepackage{dsfont}
\usepackage{amsthm}
\usepackage[usenames,dvipsnames,svgnames,table]{xcolor}
\usepackage{url}
\usepackage{nameref}

\makeatletter
\let\orgdescriptionlabel\descriptionlabel
\renewcommand*{\descriptionlabel}[1]{%
  \let\orglabel\label
  \let\label\@gobble
  \edef\@currentlabel{#1}%
  \let\label\orglabel
  \orgdescriptionlabel{#1}%
}
\makeatother

\newtheorem{theorem}{Theorem}[section]

\newtheorem{lemma}[theorem]{Lemma}

\theoremstyle{definition}

\newtheorem{algorithm}[theorem]{Algorithm}

\theoremstyle{remark}

\numberwithin{equation}{section}

\newcommand{\argmin}{\operatornamewithlimits{arg\, min}}
\newcommand{\bs}{\boldsymbol}
\newcommand{\sr}{\dot{\gamma}}
\newcommand{\sg}{\mathcal{D}}

\renewcommand{\div}{\operatorname{div}}
\newcommand{\Div}{\operatorname{Div}}
\newcommand{\prox}{\operatorname{prox}}

\newcommand{\bi}{\mathrm{Bi}}

\usepackage{enumerate}

\begin{document}

\begin{frontmatter}

\title{An Accelerated Dual Proximal Gradient Method for Applications in Viscoplasticity}
\author[uc]{Timm Treskatis\corref{corauth}\fnref{support}}
\cortext[corauth]{Corresponding author}
\fntext[support]{Supported by a UC Doctoral Scholarship and the OptALI Exchange Programme}
\ead{timm.treskatis@pg.canterbury.ac.nz}
\author[uc]{Miguel A. Moyers-Gonz\'alez}
\author[uc]{Chris J. Price}
\address[uc]{School of Mathematics and Statistics, University of Canterbury, Private Bag 4800, Christchurch 8140, New Zealand}
\journal{arXiv}

\begin{abstract}
We present a very simple and fast algorithm for the numerical solution of viscoplastic flow problems without prior regularisation. Compared to the widespread alternating direction method of multipliers (ADMM / ALG2), the new method features three key advantages: firstly, it accelerates the worst-case convergence rate from  $O(1/\sqrt{k})$ to $O(1/k)$, where $k$ is the iteration counter. Secondly, even for nonlinear constitutive models like those of Casson or Herschel-Bulkley, no nonlinear systems of equations have to be solved in the subproblems of the algorithm. Thirdly, there is no need to augment the Lagrangian, which eliminates the difficulty of choosing a penalty parameter heuristically.

In this paper, we transform the usual velocity-based formulation of viscoplastic flow problems to a dual formulation in terms of the stress. For the numerical solution of this dual problem we apply FISTA, an accelerated first-order optimisation algorithm from the class of so-called proximal gradient methods. Finally, we conduct a series of numerical experiments, focussing on stationary flow in two-dimensional square cavities.

Our results confirm that Algorithm FISTA*, the new dual-based FISTA, outperforms state-of-the-art algorithms such as ADMM / ALG2 by several orders of magnitude. We demonstrate how this speedup can be exploited to identify the free boundary between yielded and unyielded regions with previously unknown accuracy. Since the accelerated algorithm relies solely on Stokes-type subproblems and nonlinear function evaluations, existing code based on augmented Lagrangians would require only few minor adaptations to obtain an implementation of FISTA*.
\end{abstract}

\end{frontmatter}

\section{Introduction}

Viscoplasticity is a wide-spread phenomenon in both natural and man-made applications. The rich rheology of viscoplastic fluids is encountered in geophysics, considering the examples of lava flows or lahars \cite{Balmforth2000,Manville1998}. Certain types of mineral oils, mud or slurry suspensions also exhibit viscoplastic features. In the consumer goods industry, toothpaste, hair gel, tomato sauce or dough serve as classical examples of such fluids \cite{Bird1983}.

The characteristic feature of a viscoplastic fluid is its ability to resist stress in the material up to a critical threshold, the so-called yield stress $\tau_0$. This behaviour is generally due to friction-type interactions between the molecules or particles of the fluid. Consequently, viscoplastic fluids behave like a rigid material at small stress. They only start shearing like a viscous liquid if the stress exceeds the threshold posed by the yield stress.

\subsection{Mathematical Models for Viscoplastic Fluid Flows}

We consider the problem of steady, creeping viscoplastic flow in a cavity, represented by the bounded domain $\Omega \subset \mathds{R}^d$ with (Lipschitz) boundary $\partial \Omega = \Gamma$. In practice, $d\in \lbrace 2,3 \rbrace$. Our objective is to solve for functions $\bs u: \Omega \to \mathds{R}^d$, $p: \Omega \to \mathds{R}$ and $\bs \tau: \Omega \to \mathds{R}^{d\times d}_\mathrm{sym}$, representing the flow velocity, pressure and deviatoric part of the stress, respectively. Furthermore, with the symmetric gradient operator $\sg := (\nabla + \nabla^\top)/2$, we denote the strain-rate tensor by $\sr: \Omega \to \mathds{R}^{d\times d}_\mathrm{sym}$, which is linked to the flow velocity through the relation $\sr = \sg \bs u$.

The most common mathematical descriptions of viscoplastic behaviour are given by the Bingham \cite{Bingham1922}, the Casson \cite{Casson1959} and the shear-thinning Herschel-Bulkley model \cite{Herschel1926}. With viscosity or consistency parameters $\mu,\kappa > 0$ and an exponent  $1<r<2$, they can be formulated as
\begin{subequations} \label{eq:constitutive-dim}
\begin{align}
\lvert \bs \tau \rvert &\leq \tau_0 && \text{if } \sr = 0\\
\bs \tau &=
\begin{dcases*}
2 \mu \sr + \tau_0 \frac{\sr}{\lvert \sr \rvert} & (Bingham)\\
\left(\sqrt{2 \mu \lvert \sr \rvert} + \sqrt{\tau_0}\right)^2 \frac{\sr}{\lvert \sr \rvert} & (Casson)\\
2^{r-1} \kappa \lvert\sr\rvert^{r-2} \sr + \tau_0 \frac{\sr}{\lvert \sr \rvert} & (Herschel-Bulkley)
\end{dcases*}
&& \text{if }\sr \neq 0.
\end{align}
\end{subequations}
Here, $\vert \cdot \vert$ denotes the Frobenius norm on $\mathds{R}^{d \times d}_\mathrm{sym}$.

In what follows, we consider a non-dimensionalised formulation that has been re-scaled with respect to a characteristic length $L$ and velocity $U$, which reduces the dimensional constitutive relations \eqref{eq:constitutive-dim} to
\begin{subequations}\label{eq:constitutive}
\begin{align}
\lvert \bs \tau \rvert & \leq \bi && \text{if } \sr = 0\\
\bs \tau &=
\begin{dcases*}
2 \sr + \bi \frac{\sr}{\lvert \sr \rvert} & (Bingham)\\
\left(\sqrt{2 \lvert \sr \rvert} + \sqrt{\bi}\right)^2 \frac{\sr}{\lvert \sr \rvert} & (Casson)\\
2^{r-1} \lvert\sr\rvert^{r-2} \sr + \bi \frac{\sr}{\lvert \sr \rvert} & (Herschel-Bulkley)
\end{dcases*}
&& \text{if }\sr \neq 0.
\end{align}
\end{subequations}
The Bingham number $\bi := \frac{\tau_0 L}{\mu U}$ (Bingham, Casson) or $\bi := \frac{\tau_0 L^{r-1}}{\kappa U^{r-1}}$ (Herschel-Bulkley) quantifies the deviation of the viscoplastic flow from (generalised) Newtonian behaviour.

Any of these constitutive relations, along with equations for conservation of momentum and mass, yield a system for the unknown flow variables. Denoting by $\boldsymbol f: \Omega \to \mathds{R}^d$ a non-dimensionalised density of body forces, we have
\begin{align}
-\Div {\bs \tau} + \nabla p &= \boldsymbol f & &\mbox{in } \Omega \label{eq:momentum}\\
\div {\bs u} &= 0 & &\mbox{in }\Omega. \label{eq:incompressibility}
\intertext{To close the system, we incorporate the boundary condition}
\bs u &= \bs u_\mathrm{D} & &\mbox{on } \Gamma \label{eq:bcon},
\end{align}
where $\boldsymbol u_\mathrm{D}: \Gamma_\mathrm{D} \to \mathds{R}^d$ is given. We use the notation $\div$ (resp. $\Div$) for the (rowwise) divergence operator.

\subsection{Variational Formulation}

In the following, we use boldface letters for spaces to denote $d$-fold Cartesian products, e.g. for a space $A$ we write $\bs A := A^d$. To obtain a mathematically rigorous formulation of the viscoplastic flow problem \eqref{eq:constitutive}-\eqref{eq:bcon} in Sobolev spaces, we consider
\begin{align*}
U &:= \bs W^{1,r}(\Omega)\\
U_{0*} &:= \Set{\bs u \in \bs W^{1,r}(\Omega) | \div \bs u = 0}\\
U_{*0} &:= \Set{\bs u \in \bs W^{1,r}(\Omega) | \bs u\rvert_{\Gamma} = 0}\\
U_{00} &:= \Set{\bs u \in \bs W^{1,r}(\Omega) | \div \bs u = 0 \text{ and } \bs u\rvert_{\Gamma} = 0},
\end{align*}
with $r=2$ in the Bingham and Casson settings. We use the dual of the latter space to fix the inhomogeneity $\bs f \in U_{00}^*$ and we pick boundary values $\bs u_\mathrm{D} \in U_\mathrm{D}$, where
\begin{equation*}
U_\mathrm{D} := \Set{ \bs u_\mathrm{D} \in \bs W^{1-1/r,r}(\Gamma) | \int\limits_\Gamma \bs{u}_\mathrm{D} \cdot \bs n \, \mathrm{d}s = 0 }.
\end{equation*}
Furthermore, we define the convex set of admissible solutions
\begin{equation*}
U_{0\mathrm{D}} := \Set{\bs u \in \bs W^{1,r}(\Omega) | \div \bs u = 0 \text{ and } \bs u\rvert_\Gamma = \bs u_\mathrm{D} }.
\end{equation*}
For the strain-rate and stress tensors, we will also need spaces of symmetric matrices whose entries satisfy an integrability condition of order $r$ or $r^*$, respectively, where $1/r + 1/r^* = 1$:
\begin{equation*}
Q := L^r (\Omega)^{d\times d}_{\mathrm{sym}} \qquad \qquad S := L^{r^*} (\Omega)^{d \times d}_\mathrm{sym}.
\end{equation*}

By generalising the ideas of Duvaut and Lions \cite{Duvaut1972,Duvaut1976} and Huilgol and You \cite{Huilgol2005}, we conclude that the system \eqref{eq:constitutive}-\eqref{eq:bcon} is a strong formulation of the following variational inequality problem of the second kind: find $\bs u \in  U_{0\mathrm{D}}$ such that for all test velocity fields $\bs v \in U_{0\mathrm{D}}$
\begin{equation} \label{eq:vi}
a(\sg \bs u , \sg \bs v - \sg \bs u) + j(\sg \bs v) - j(\sg \bs u) \geq \langle \bs f, \bs v - \bs u \rangle_{U_{00}^*,U_{00}}.
\end{equation}
This variational inequality is composed of the elliptic form $a: Q\times Q \to \mathds{R}$,
\begin{align*}
a(\bs\sr, \bs{\dot{\delta}}) &:= 2\int\limits_\Omega \bs\sr : \bs{\dot{\delta}} \, \mathrm{d}x & & \text{(Bingham)}\\
a(\bs\sr, \bs{\dot{\delta}}) &:= 2\int\limits_\Omega \bs\sr : \bs{\dot{\delta}} \, \mathrm{d}x + 2 \sqrt{2\bi} \int\limits_\Omega \frac{\bs\sr}{\sqrt{\lvert \bs\sr \rvert}} : \bs{\dot{\delta}} \, \mathrm{d}x & & \text{(Casson)}\\
a(\bs\sr, \bs{\dot{\delta}}) &:= 2^{r-1} \int\limits_\Omega \lvert\bs\sr\rvert^{r-2} \bs\sr : \bs{\dot{\delta}} \, \mathrm{d}x & & \text{(Herschel-Bulkley)}
\end{align*}
the nonsmooth functional $j: Q \to \mathds{R}$,
\begin{equation*}
j(\bs\sr) := \bi \int\limits_\Omega \lvert \bs\sr \rvert \,\mathrm{d}x
\end{equation*}
and, on the right-hand side, a duality pairing between $U_{00}$ and its dual, which can be represented as
\begin{equation*}
\langle \bs f, \bs v - \bs u \rangle_{U_{00}^*,U_{00}} = \int\limits_\Omega \bs f \cdot (\bs v - \bs u)\, \mathrm{d}x
\end{equation*}
provided that $\bs f \in \bs L^{r^*}(\Omega)$. The colon represents the Frobenius inner product of two $d \times d$ matrices, the dot the scalar product of two vectors in $\mathds{R}^d$.

It is an important observation that for each of the three viscoplastic models, the term $a(\sg \bs u , \sg \bs v - \sg \bs u)$ possesses special structure: with the functional $b : Q \to \mathds{R}$ defined by
\begin{align*}
b(\bs\sr) &:= \int\limits_\Omega \lvert \bs\sr \rvert^2 \, \mathrm{d}x & & \text{(Bingham)}\\
b(\bs\sr) &:= \int\limits_\Omega \lvert \bs\sr \rvert^2 \, \mathrm{d}x + \frac{4\sqrt{2\bi}}{3} \int\limits_\Omega \lvert \bs\sr \rvert^{3/2} \, \mathrm{d}x & & \text{(Casson)}\\
b(\bs\sr) &:= \frac{2^{r-1}}{r} \int\limits_\Omega \lvert \bs\sr \rvert^r \, \mathrm{d}x & & \text{(Herschel-Bulkley)}
\end{align*}
we may write 
\begin{equation*}
a(\sg \bs u , \sg \bs v - \sg \bs u) = \left\langle \nabla_{\bs u} b (\sg \bs u), \sg (\bs v - \bs u) \right\rangle_{Q^*,Q}
\end{equation*}
as a directional derivative of $b\circ \sg$ at $\bs u$ in direction $\bs v - \bs u$. Consequently, we may identify the variational inequality \eqref{eq:vi} as a first-order optimality condition of the convex, and hence equivalent minimisation problem
\begin{equation} \label{eq:min} \tag{VP}
\min_{\bs u \in U} b(\sg \bs u) + j(\sg \bs u) - \langle \bs f, \bs u \rangle_{\bs L^{r^*}(\Omega),\bs L^r(\Omega)} + \iota_{U_{0\mathrm{D}}}(\bs u),
\end{equation}
where the indicator functional
\begin{equation*}
\iota_{U_{0\mathrm{D}}}(\bs u) =
\begin{dcases*}
0 & if $\bs u \in U_{0\mathrm{D}}$\\
+\infty & if $\bs u \notin U_{0\mathrm{D}}$
\end{dcases*}
\end{equation*}
enforces the incompressibility constraint \eqref{eq:incompressibility} and the Dirichlet boundary condition \eqref{eq:bcon}.

For full details of the derivation of Problem \eqref{eq:min}, and results regarding existence and uniqueness of solutions, we refer to \cite[Ch~4]{Treskatis2016}.

While the Bingham and Casson flow problems are posed in Hilbert spaces, the Herschel-Bulkley model demands for a mathematical treatment in more general Banach spaces. Despite extra theoretical challenges, a very practical consequence of this fact is that a numerical optimisation algorithm would include a series of nonlinear subproblems since the exponents $r$ and $r^*$ are different and no longer equal to two. To circumvent these difficulties, it has become common practice in viscoplasticity (see e.g. the appendix of \cite{Huilgol2005} and the references therein) to discretise the problem first, and then solve a discretised approximation of \eqref{eq:min}. This way, all spaces $U$, $Q$, $S$ are replaced by finite-dimensional spaces $U_h$, $Q_h$, $S_h$, e.g. spaces of finite elements, for which the Hilbert-space structure can be recovered.

Strictly speaking, we would now have to differentiate between such finite-dimensional spaces, whenever we refer to the Herschel-Bulkley problem, while we concurrently work with the original function spaces for the Bingham and Casson problems. In an attempt to provide a clearer picture, we will not show this distinction in our notation and hide the subscript $h$ for now. We shall however emphasise that in order to re-establish a mathematically rigorous formulation for Herschel-Bulkley fluids, one has to replace all variables with their finite-dimensional counterparts.

\subsection{State-of-the-Art Approaches to Solving (\ref{eq:min})}

The community of numerical analysts in viscoplasticity can essentially be divided into two groups: some authors approximate the constitutive relations \eqref{eq:constitutive} with a more regular formulation, while others leave the genuinely nonsmooth nature of the original problem \eqref{eq:min} unaltered.

For formulations of the first kind, we exemplarily mention the Bercovier-Engelman model \cite{Bercovier1980}, the Papanastasiou regularisation \cite{Papanastasiou1987}, bi-viscosity formulations of Tanner and Milthorpe \cite{Tanner1983} and De los Reyes and Gonz{\'a}lez Andrade \cite{Gonzalez-Andrade2008,Reyes2009,Reyes2010,Reyes2012} and the penalty approach of Glowinski and collaborators \cite{Dean2007,Glowinski2003,Dacorogna2004}. It can be seen as the main advantage of these approximations that they allow for very efficient numerical methods of Newton-type, with a fast, locally superlinear or quadratic convergence rate. However, for many practical applications it can be problematic that such approximate solutions do not generally reflect the characteristic features of the exact solution. For instance, Moyers-Gonz{\'a}lez and Frigaard \cite{Moyers2004} demonstrate that even under arbitrarily small excitations, solutions to regularised models predict slow flow although the actual flow rate is exactly zero, see also \cite{Frigaard2005,Balmforth2014}. Another well-known issue arises from the difficulties of recovering the yielded and unyielded flow regions under a regularised problem formulation. In fact, it appears that convergence of these approximate plug and shear regions has not been proved yet, and such a statement may not even hold \cite{Balmforth2014}.

In our work, we employ the genuinely nonsmooth formulation \eqref{eq:min}. Classical algorithms for the numerical solution of this convex optimisation problem stem from the framework of augmented Lagrangians. A set of augmented Lagrangian methods ALG1--ALG4 was originally proposed by Fortin and Glowinski \cite{Fortin1983}. They are based on the idea to introduce $\bs \sr = \sg \bs u$ as a constraint to the minimisation problem, the violation of which is then penalised with an extra quadratic term:
\begin{equation} \tag{$\text{VP}^\prime_\varrho$}\label{eq:min-split-aug}
\min_{(\bs u,\bs \sr) \in U\times Q} b(\bs \sr) + j(\bs \sr) - \langle \bs f, \bs u \rangle_{\bs L^{r^*}(\Omega),\bs L^r(\Omega)} + \iota_{U_{0\mathrm{D}}}(\bs u) + \frac{\varrho}{2}\left\lVert \bs \sr - \sg \bs u \right\rVert_Q^2 \qquad \text{subject to } \bs \sr = \sg \bs u.
\end{equation}
The problems \eqref{eq:min} and \eqref{eq:min-split-aug} are clearly equivalent for any $\varrho \geq 0$.

As another equivalent alternative to \eqref{eq:min-split-aug}, one frequently encounters the formulation
\begin{equation*}
\min_{(\bs u,\bs \sr) \in U\times Q} b(\sg \bs u) + j(\bs \sr) - \langle \bs f, \bs u \rangle_{\bs L^{r^*}(\Omega),\bs L^r(\Omega)} + \iota_{U_{0\mathrm{D}}}(\bs u) + \frac{\varrho}{2}\left\lVert \bs \sr - \sg \bs u \right\rVert_Q^2 \qquad \text{subject to } \bs \sr = \sg \bs u
\end{equation*}
in the literature on viscoplastic flow, or the corresponding first-order optimality conditions. We refer to \cite{Duvaut1972,Duvaut1976,Reyes2010,Dean2007} for a more detailed discussion in the Bingham setting. This approach has the slight disadvantage that it leads to a rather complicated dual problem. As shown, e.g. in \cite{Reyes2010}, it has the form of an elliptic optimal control problem with pointwise inequality constraint on the control. The dual multiplier has the dimension of a stress and can be interpreted as a plastic contribution to the extra stress tensor. We will consider the formulation in \eqref{eq:min-split-aug}, where the Lagrange multiplier can be identified with the stress $\bs\tau$ itself.

ALG1, the generalised Uzawa method \cite{Uzawa1958} applied to the augmented Lagrangian that corresponds to \eqref{eq:min-split-aug}, attempts to find a saddle point by minimising in the variables $(\bs u,\bs \sr)$ and then taking a step along the dual gradient in order to maximise with respect to the Lagrange multiplier $\bs \tau$. Since minimising jointly in $\bs u$ and $\bs \sr$ is as difficult as solving the original problem, every iteration of ALG1 is very costly.

In ALG2, the alternating direction method of multipliers (ADMM), the exact solution of each subproblem is waived in favour of a minimisation in $\bs u$ only and a subsequent minimisation in $\bs \sr$ only. ALG2 is related to the Douglas-Rachford splitting algorithm. For further details, we refer to Glowinski's recent review \cite{Glowinski2014}, the references therein as well as \cite{Gabay1976}. Of all augmented Lagrangian methods, ALG2 is by far the most popular one and can be considered as the standard genuinely nonsmooth approach to simulating viscoplastic fluid flows. Even though the convergence analysis of ALG2 is difficult, it is meanwhile well-known that even under assumptions that are too strong for Problem \eqref{eq:min}, only a sublinear convergence rate of $O(1/k)$ in the dual objective functional can be guaranteed \cite[Thm~1]{Goldstein2014}. This corresponds to a convergence rate of only $O(1/\sqrt{k})$ for the primal iterates, which represent the velocity and strain rate.

ALG3 is a counterpart of the Peaceman-Rachford splitting method. Also known as alternating minimisation algorithm (AMA), it expands on ALG2 by adding an update of the dual variable $\bs \tau$ not just after the minimisation in $\bs \sr$, but also after the minimisation in $\bs u$. ALG4 was designed as an analogue of the $\theta$-method \cite{Fortin1983,Glowinski1989}, but appears to be rarely applied.

Convex optimisation problems of a very similar structure to \eqref{eq:min} are also encountered in a variety of other disciplines, such as signal and image processing \cite{Ngwa2003,Frigaard2006}, machine learning, statistics or mathematical finance. The field of nonsmooth convex optimisation consequently receives considerable attention and many alternative numerical methods have been derived in past years. To us, the class of so-called proximal gradient algorithms appears particularly well-suited for the solution of Problem \eqref{eq:min}, as these methods can readily exploit the structure of composite convex objectives: terms which are smooth in the sense that they possess a Lipschitz-continuous gradient are linearised, while nonsmooth terms are left unchanged.

In recent years, interest in proximal algorithms has reached an unprecedented extent with the re-discovery and further development of fast or accelerated methods. Building upon ideas of Nesterov \cite{Nesterov1983,Nesterov2007}, Beck and Teboulle \citep{Beck2009} constructed a fast iterative shrinkage-thresholding algorithm (FISTA), an efficient fast variant of the basic proximal gradient method. At every iteration, FISTA extrapolates from the current and immediate past iterate to accelerate the convergence of the method. While the conventional, unaccelerated method decreases the value of the objective at a rate of order $O(1/k)$, where $k$ is the iteration counter, FISTA achieves $O(1/k^2)$ at negligible extra computational cost. In particular, it requires no accumulation of the iteration history, apart from the last iterate.

Very recently, Beck and Teboulle considered abstract composite problems of the form
\begin{equation*}
\min_{\bs u} f(\bs u) + g(\sg \bs u)
\end{equation*}
with a strongly convex functional $f$ and a general convex functional $g$. In \cite{Beck2014} (see also \cite{Beck2009a}), they apply FISTA to the dual problem to achieve $O(1/k)$-convergence of the primal sequence, in contrast to $O(1/\sqrt{k})$ without acceleration. Given that problem \eqref{eq:min} is of the similar form
\begin{equation*}
\min_{\bs u} f(\sg \bs u) + g(\bs u)
\end{equation*}
again with a strongly convex $f$ and convex $g$, our approach exhibits some analogies to \cite{Beck2014}. The similarity between our and their problem allows us to apply some of their ideas in our work, in particular the concept of applying FISTA to the dual problem. However, there are also two important differences: firstly, the different structure of the problem gives rise to some additional difficulties and requires further technical assumptions in order to guarantee the applicability and convergence of FISTA. Secondly, Beck and Teboulle restrict themselves to problems in $\mathds{R}^n$, while we consider a general function space setting.

Even though our approach is applicable to very general composite convex problems, we will only focus on the application of the methodology to Problem \eqref{eq:min}. The more abstract mathematical details can be found in \cite{Treskatis2016}.

We shall point out that acceleration techniques have also been studied for the augmented Lagrangian methods ALG2 and ALG3, see \cite{Goldstein2014,Goldfarb2013}. However, quantitative rates of convergence are very difficult to prove for these methods. At this stage, it appears that there are only isolated results available that rely on fairly restrictive assumptions or unconventional measures of convergence \cite{Goldstein2014, He2012, He2012a, Deng2012, Hong2012}. 

For an extensive review of state-of-the-art gradient-based optimisation methods for abstract composite convex problems, we refer to \cite{Burger2014}.

\paragraph*{Outline}

In Section \ref{sec:dual}, we introduce a novel dual formulation of Problem \eqref{eq:min}. For its numerical solution, we present the accelerated dual proximal gradient method FISTA* in Section \ref{sec:algorithm} along with some key properties. Finally, we conclude with numerical results and a discussion in Section \ref{sec:results}.

\section{Dual Formulation}
\label{sec:dual}

\subsection{Derivation of the Dual Problem}

We consider Problem \eqref{eq:min-split}, the split formulation of the viscoplastic flow problem with no penalty term added:
\begin{equation} \tag{$\text{VP}^\prime_0$}\label{eq:min-split}
\min_{(\bs u,\bs \sr) \in U\times Q} b(\bs \sr) + j(\bs \sr) - \langle \bs f, \bs u \rangle_{\bs L^{r^*}(\Omega),\bs L^r(\Omega)} + \iota_{U_{0\mathrm{D}}}(\bs u) \qquad \text{subject to } \bs \sr = \sg \bs u.
\end{equation}
By introducing a Lagrange multiplier (dual variable), which has a physical interpretation as an admissible stress tensor $\bs \tau$ \cite{Huilgol2005}, we can re-write this problem equivalently as
\begin{equation*}
\max_{\bs\tau \in S} \min_{(\bs u,\bs \sr) \in U\times Q} b(\bs \sr) + j(\bs \sr) - \langle \bs f, \bs u \rangle_{\bs L^{r^*}(\Omega),\bs L^r(\Omega)} + \iota_{U_{0\mathrm{D}}}(\bs u) - \left\langle \bs\tau, \bs \sr - \sg \bs u\right\rangle_{Q*,Q}.
\end{equation*}
We now solve the inner minimisation problem to eliminate the primal variables $\bs u$ and $\bs \sr$. Re-arranging yields
\begin{equation}\label{eq:dual-aux}
\max_{\bs\tau \in S} \left(\min_{\bs \sr \in Q} \Set{b(\bs \sr) + j(\bs \sr) - \left\langle \bs\tau, \bs \sr\right\rangle_{Q*,Q} } + \min_{\bs u \in U} \Set{- \langle \bs f, \bs u \rangle_{\bs L^{r^*}(\Omega),\bs L^r(\Omega)} + \iota_{U_{0\mathrm{D}}}(\bs u) + \left\langle \bs\tau, \sg \bs u\right\rangle_{Q*,Q}} \right).
\end{equation}

For the second minimisation problem in \eqref{eq:dual-aux} we infer
\begin{equation*}
\min_{\bs u \in U} \Set{- \langle \bs f, \bs u \rangle_{\bs L^{r^*}(\Omega),\bs L^r(\Omega)} + \iota_{U_{0\mathrm{D}}}(\bs u) + \left\langle \bs\tau, \sg \bs u\right\rangle_{Q*,Q}}
=
\begin{dcases}
0 & \text{if } \bs\tau \in C\\
- \infty & \text{otherwise}.
\end{dcases}
= - \iota_C(\bs\tau),
\end{equation*}
with the set
\begin{equation*}
C := \Set{\bs\tau\in S | \left\langle \bs\tau, \sg \bs u\right\rangle_{Q*,Q} - \langle \bs f, \bs u \rangle_{\bs L^{r^*}(\Omega),\bs L^r(\Omega)} = 0 \text{ for all } \bs u \in U_{00}}.
\end{equation*}
We point out that this condition on $\bs\tau$ is a very compact variational formulation of \eqref{eq:momentum}. This becomes more obvious if we introduce the pressure $p$ as a Lagrange multiplier for the constraint $\div \bs u = 0$, which is implicitly contained in the definition of the space $U_{00}$. Then the set $C$ contains all $\bs \tau \in S$, for which there exist $p \in L^{r^*} (\Omega)$ such that
\begin{equation*}
\left\langle \bs\tau, \sg \bs u\right\rangle_{Q*,Q} - \langle p, \div \bs u \rangle_{L^{r^*}(\Omega),L^r(\Omega)} = \langle \bs f, \bs u \rangle_{\bs L^{r^*}(\Omega),\bs L^r(\Omega)}
\end{equation*}
for all test functions $\bs u \in U_{*0}$.

Let us now turn to the first minimisation problem in \eqref{eq:dual-aux}. Separate calculations for the three different constitutive models reveal the following explicit result \cite[Ch~6]{Treskatis2016}:
\begin{equation*}
\min_{\bs \sr \in Q} \Set{b(\bs \sr) + j(\bs \sr) - \left\langle \bs\tau, \bs \sr\right\rangle_{Q*,Q} } = - F(\bs\tau)
\end{equation*}
with
\begin{equation}
F(\bs\tau) = 
\begin{dcases*}
\frac{1}{4}\int\limits_\Omega \left(\lvert\bs\tau\rvert - \bi\right)_+^2\,\mathrm{d}x & (Bingham)\\
\frac{1}{4}\int\limits_\Omega \left(\sqrt{\lvert\bs\tau\rvert} - \sqrt{\bi}\right)_+^3 \left( \sqrt{\lvert\bs\tau\rvert} + \frac{1}{3}\sqrt{\bi} \right)\,\mathrm{d}x & (Casson)\\
\frac{1}{2r^*}\int\limits_\Omega \left(\lvert\bs\tau\rvert - \bi\right)_+^{r^*}\,\mathrm{d}x & (Herschel-Bulkley).
\end{dcases*}
\end{equation}
Here, we employed the notation $(\cdot)_+ = \max\Set{0,\cdot}$.

In summary, \eqref{eq:dual-aux} leads to the following dual formulation of the viscoplastic flow problem:
\begin{equation}\tag{$\text{VP}^*$}\label{eq:dual}
\min_{\bs \tau \in S} F(\bs\tau) + \iota_C(\bs\tau)
\end{equation}
or simply
\begin{equation*}
\min_{\bs \tau \in C} F(\bs\tau).
\end{equation*}

\subsection{Properties of the Dual Problem}

The dual viscoplastic flow problem \eqref{eq:dual} is a composite convex optimisation problem, where both terms possess special properties.

For all three constitutive models, the functional $F$ is (Fr{\'e}chet-)differentiable. After short calculations we find that its gradient with respect to $\bs\tau$ is given by
\begin{equation}\label{eq:dual-grad}
\nabla F(\bs\tau) = 
\begin{dcases*}
\frac{1}{2}\left(\lvert\bs\tau\rvert - \bi\right)_+\frac{\bs\tau}{\lvert\bs\tau\rvert} & (Bingham)\\
\frac{1}{2}\left(\sqrt{\lvert\bs\tau\rvert} - \sqrt{\bi}\right)_+^2\frac{\bs\tau}{\lvert\bs\tau\rvert} & (Casson)\\
\frac{1}{2}\left(\lvert\bs\tau\rvert - \bi\right)_+^{r^*-1}\frac{\bs\tau}{\lvert\bs\tau\rvert} & (Herschel-Bulkley).
\end{dcases*}
\end{equation}
In fact, it holds that $\nabla F (\bs\tau) = \bs\sr$: by re-arranging this equation for the stress $\bs\tau$, one recovers the constitutive models of \eqref{eq:constitutive}.

Additionally, for Bingham and Casson fluids, $\nabla F$ is globally Lipschitz-continuous, where the smallest possible Lipschitz constant is $L = 1/2$. This could either be verified directly, but it follows more easily from classical relations between primal and dual problems in convex optimisation \cite[Lemma~3.6]{Treskatis2016}.

Note that since $1<r<2$ for shear-thinning Herschel-Bulkley fluids, the dual exponent $r^* > 2$. Therefore, in this case $\nabla F$ is only Lipschitz-continuous on bounded subsets of the stress space $S$. However, as long as an algorithm for the solution of Problem \eqref{eq:dual} does not diverge, all stress iterates will remain bounded and we can still find a Lipschitz constant $L>0$ for the set of all iterates. In contrast to the Bingham and Casson settings, it appears impossible in practice to determine this constant a priori. Instead, for Herschel-Bulkley flow problems we will have to compute estimates for $L$ numerically.

The second term in \eqref{eq:dual} is an indicator functional, which imposes the momentum equation in form of the constraint $\bs \tau \in C$. In convex optimisation, indicator functionals of convex sets serve as a prototypical example for \emph{functionals of simple structure}. A functional $G$ is said to be simple if its so-called Moreau proximal map \cite{Moreau1965},
\begin{equation} \label{eq:proximal}
\prox_G (\bs\tau) := \argmin_{\bs \sigma \in S} \Set{ G(\bs \sigma) + \frac{1}{2} \lVert \bs\sigma - \bs\tau \rVert_S^2 },
\end{equation}
can be evaluated in a computationally efficient manner. If $G$ is an indicator function of the set $C$, then its proximal map equates to a projection of $\bs \tau$ onto $C$.

We are now in the position to proceed with the numerical solution of \eqref{eq:dual}.

\section{Accelerated Dual Proximal Gradient Method}
\label{sec:algorithm}

Our objective for this section is to present algorithmic approaches to solving the primal viscoplastic flow problem \eqref{eq:min} through its dual \eqref{eq:dual}. First, we briefly review the fundamentals of FISTA and then apply this method to \eqref{eq:dual}, in order to derive the basic algorithm FISTA* for the solution of \eqref{eq:min}. After a comparison with the alternating direction method of multipliers (ADMM / ALG2), we present our results on convergence rates. Finally, we comment on some aspects of the implementation.

\subsection{Proximal Gradient Methods}

The fast iterative shrinkage-thresholding algorithm (FISTA) belongs to the class of proximal gradient methods. Its name stems from the fact that for certain convex functionals $G$, the proximal map \eqref{eq:proximal} becomes a shrinkage-thresholding operator. However, it has meanwhile become common practice \cite{Beck2009} to slightly misuse the term \emph{iterative shrinkage-thresholding algorithm} even if the proximal map of $G$ is of different form, as is the case when $G = \iota_C$ actually yields an iterative projection algorithm.

The most basic proximal gradient method, or ISTA, for the solution of the composite convex problem
\begin{equation*}
\min_{\bs \tau \in S} F(\bs\tau) + G(\bs\tau),
\end{equation*}
where $F$ has a Lipschitz-continuous gradient with Lipschitz constant $L$, reads as follows \cite{Beck2009}:

\begin{algorithm}[ISTA / Proximal Gradient Method]\label{alg:ista}
\emph{Input:} $\bs\tau^{(0)}\in S$\\
\emph{Initialisation:} $k = 1$
\begin{description}
\item[(ISTA.1)\label{it:ista1}] Set $L^{(k)} = L$ or find a valid $L^{(k)} > 0$ by backtracking and evaluate
\begin{equation*}
\bs\tau^{(k)} = \prox_{\frac{1}{L^{(k)}}G} \left( \bs\tau^{(k-1)} - \frac{1}{L^{(k)}}\nabla F (\bs\tau^{(k-1)})\right)
\end{equation*}
\item[(ISTA.2)\label{it:ista2}] \emph{If} the algorithm has converged \emph{Then Return} $\bs\tau^{(k)}$ and \emph{Stop}.
\item[(ISTA.3)\label{it:ista3}] Set $k \leftarrow k+1$ and \emph{Go To} \ref{it:ista1}.
\end{description}
\end{algorithm}

A valid estimate $L^{(k)}$ of the Lipschitz constant $L$ means that $L^{(k)}$ and the corresponding iterate $\bs\tau^{(k)}$ must satisfy the descent criterion
\begin{equation}\label{eq:descent-ista}
F(\bs\tau^{(k)}) \leq F(\bs\tau^{(k-1)}) + \left\langle \nabla F (\bs\tau^{(k-1)}), \bs\tau^{(k)}-\bs\tau^{(k-1)} \right\rangle_{S^*,S} + \frac{L^{(k)}}{2} \left\lVert \bs\tau^{(k)}-\bs\tau^{(k-1)}\right\rVert_S^2,
\end{equation}
which holds in particular for all $L^{(k)} \geq L$ \cite[Lemma~2.1]{Beck2009}. To find a valid parameter $L^{(k)}$ in an algorithmic fashion, one would start from a trial value $\tilde{L}$ and increase this value by a certain percentage until \eqref{eq:descent-ista} is satisfied. We will present the full details below in Algorithm \ref{alg:bt}.

The following convergence result for ISTA is established in \cite{Beck2009}, Remark 2.1 and Theorem 3.1:

\begin{theorem}[Convergence Rate of ISTA]\label{thm:cvgce-ista}
Let the sequence $(\bs\tau^{(k)})_k$ be generated by ISTA, where the sequence $(L^{(k)})_k$ is non-decreasing and bounded. If $\bs{\bar{\tau}}\in\argmin_{\bs\tau\in S} \Set{F(\bs\tau)+G(\bs\tau)}$ is any solution of the minimisation problem, then
\begin{equation*}
\left( F(\bs\tau)+G(\bs\tau) \right) - \left( F(\bs{\bar{\tau}})+G(\bs{\bar{\tau}}) \right) \leq O(1/k).
\end{equation*}
\end{theorem}

He and Yuan \cite{He2012a} prove a corresponding $O(1/k)$-convergence result for the alternating direction method of multipliers (ADMM / ALG2), but it appears that their stronger assumptions on the problem compared to Theorem \ref{thm:cvgce-ista} cannot generally be relaxed.

It turns out that this convergence rate of ISTA and ADMM / ALG2 for the class of problems with smooth convex $F$ and convex $G$ is not optimal in the sense of complexity theory (cf \cite[pp~4-7]{Nesterov2013}). In fact, there exist algorithms which converge like $O(1/k^2)$, while it can be shown \cite{Nesterov1983} that no higher rate is achievable for this entire class of optimisation problems. The fast iterative shrinkage-thresholding algorithm FISTA constitutes an example of a method which is optimal in this sense. Such accelerated or inertial algorithms have only recently attracted greater attention in convex optimisation.

Compared to ISTA, FISTA comprises an additional extrapolation step. It combines the last two iterates $\bs\tau^{(k)}$ and $\bs\tau^{(k-1)}$ in just the right amount to compute a so-called leading point $\bs{\hat{\tau}}^{(k)}$. For the next iteration, the functionals $F$ and $G$ and the gradient $\nabla F$ are then evaluated at this leading point $\bs{\hat{\tau}}^{(k)}$ instead of $\bs\tau^{(k)}$ to find a shortcut to the minimum. 

\begin{algorithm}[FISTA / Accelerated Proximal Gradient Method]\label{alg:fista}
\emph{Input:} $\bs\tau^{(0)}\in S$\\
\emph{Initialisation:} $k = 1$, $t^{(1)}=1$, $\bs{\hat{\tau}}^{(1)} = \bs\tau^{(0)}$
\begin{description}
\item[(FISTA.1)\label{it:fista1}] Set $L^{(k)} = L$ or find a valid $L^{(k)} > 0$ by backtracking and evaluate
\begin{equation*}
\bs\tau^{(k)} = \prox_{\frac{1}{L^{(k)}}G} \left( \bs{\hat{\tau}}^{(k)} - \frac{1}{L^{(k)}}\nabla F (\bs{\hat{\tau}}^{(k)})\right)
\end{equation*}
\item[(FISTA.2)\label{it:fista2}] \emph{If} the algorithm has converged \emph{Then Return} $\bs\tau^{(k)}$ and \emph{Stop}.
\item[(FISTA.3)\label{it:fista3}] Compute
\begin{equation*}
t^{(k+1)} = \frac{1+\sqrt{1+4t^{(k)2}}}{2}
\end{equation*}
to update the leading point
\begin{equation*}
\bs{\hat{\tau}}^{(k+1)} = \bs\tau^{(k)} + \frac{t^{(k)}-1}{t^{(k+1)}} \left( \bs\tau^{(k)} - \bs\tau^{(k-1)} \right).
\end{equation*}
\item[(FISTA.4)\label{it:fista4}] Set $k \leftarrow k+1$ and \emph{Go To} \ref{it:fista1}.
\end{description}
\end{algorithm}

Since FISTA introduces the additional variable $\bs{\hat{\tau}}^{(k)}$, its memory footprint is larger than the one of ISTA. However, the computational cost of the extrapolation step \ref{it:fista3}, which does not appear in ISTA, is only marginal: it solely requires the evaluation of a linear combination. 

Similar to \eqref{eq:descent-ista}, if no Lipschitz constant $L$ is known, then the estimate $L^{(k)}$ must at least satisfy the descent criterion
\begin{equation}\label{eq:descent-fista}
F(\bs\tau^{(k)}) \leq F(\bs{\hat{\tau}}^{(k)}) + \left\langle \nabla F (\bs{\hat{\tau}}^{(k)}), \bs\tau^{(k)}-\bs{\hat{\tau}}^{(k)} \right\rangle_{S^*,S} + \frac{L^{(k)}}{2} \left\lVert \bs\tau^{(k)}-\bs{\hat{\tau}}^{(k)}\right\rVert_S^2.
\end{equation}
We recall that the Bingham and Casson models allow us to set $L^{(k)} \equiv L = 1/2$. For Herschel-Bulkley fluids, the following backtracking strategy \cite{Beck2009} can be employed in step \ref{it:fista1}:

\begin{algorithm}[Backtracking]\label{alg:bt}
\emph{Input:} $\bs{\hat{\tau}}^{(k)}$, a trial value for $L^{(k)}>0$, a magnifying factor $\eta > 1$
\begin{description}
\item[(BT.1)\label{it:bt1}] Evaluate
\begin{equation*}
\bs\tau^{(k)} = \prox_{\frac{1}{L^{(k)}}G} \left( \bs{\hat{\tau}}^{(k)} - \frac{1}{L^{(k)}}\nabla F (\bs{\hat{\tau}}^{(k)})\right)
\end{equation*}
\item[(BT.2)\label{it:bt2}] \emph{If} \eqref{eq:descent-fista} holds \emph{Then Return} $\bs\tau^{(k)}$ and \emph{Stop}.
\item[(BT.3)\label{it:bt3}] Set $L^{(k)} \leftarrow \eta L^{(k)}$ and \emph{Go To} \ref{it:bt1}.
\end{description}
\end{algorithm}

For both fixed and variable $L^{(k)}$, Beck and Teboulle derive the following result \cite{Beck2009} (Remark 2.1 and Theorem 4.4):

\begin{theorem}[Convergence Rate of FISTA]\label{thm:cvgce-fista}
Let the sequence $(\bs\tau^{(k)})_k$ be generated by FISTA, where the sequence $(L^{(k)})_k$ is non-decreasing and bounded. If $\bs{\bar{\tau}}\in\argmin_{\bs\tau\in S} \Set{F(\bs\tau)+G(\bs\tau)}$ is any solution of the minimisation problem, then
\begin{equation*}
\left( F(\bs\tau)+G(\bs\tau) \right) - \left( F(\bs{\bar{\tau}})+G(\bs{\bar{\tau}}) \right) \leq O(1/k^2).
\end{equation*}
\end{theorem}

\subsection{The Accelerated Dual Proximal Gradient Method FISTA*}

Let us turn towards the application of FISTA to the solution of viscoplastic flow problems. We will now derive how the proximal map in step \ref{it:fista1} can be evaluated in practice and how the primal variables $\bs u$ and $\bs\sr$, representing the flow velocity and strain rate, can be recovered.

\begin{lemma} \label{lem:prox}
The assignment $\bs\tau^{(k)} = \prox_{\frac{1}{L}G} (\bs{\hat{\tau}}^{(k)} - \nabla F (\bs{\hat{\tau}}^{(k)})/L)$ is equivalent to the operations
\begin{align}
\bs{\hat{\sr}}^{(k)} &= \argmin_{\bs \sr \in Q} \Set{ b(\bs\sr) + j(\bs\sr) - \left\langle\bs{\hat{\tau}}^{(k)},\bs\sr\right\rangle_{Q*,Q} } \label{eq:sqce-sr}\\
\bs{\hat{u}}^{(k)} &= \argmin_{\bs u \in U} \Set{ -\langle \bs f, \bs u\rangle_{\bs L^{r^*}(\Omega),\bs L^r(\Omega)} + \iota_{U_{0\mathrm{D}}}(\bs u) + \frac{1}{2L} \left\lVert \sg \bs u - \left(\bs{\hat{\sr}}^{(k)} - L \bs{\hat{\tau}}^{(k)}\right) \right\rVert_Q^2 } \label{eq:sqce-u}\\
\bs{\tau}^{(k)} &= \hat{\bs{\tau}}^{(k)} + \frac{1}{L} \left( \sg \bs {\hat{u}}^{(k)} - \bs{\hat{\sr}}^{(k)} \right) \label{eq:sqce-tau}
\end{align}
\begin{proof}
We refer to \cite[Ch~3]{Treskatis2016}.
\end{proof}
\end{lemma}

In \eqref{eq:sqce-sr} and \eqref{eq:sqce-u}, the quantities $\bs{\hat{\sr}}^{(k)}$ and $\bs{\hat{u}}^{(k)}$ represent a strain rate and a velocity, respectively. They are, however, evaluated based on the extrapolated leading point $\bs{\hat{\tau}}^{(k)}$, not the actual stress iterate $\bs{\tau}^{(k)}$. The convergence result for FISTA is valid for this sequence $(\bs{\tau}^{(k)})_k$ only, while such a convergence result may not hold for the auxiliary sequence $\bs{\hat{\tau}}^{(k)}$.

We suggest three definitions for a primal sequence $(\bs u^{(k)},\bs \sr^{(k)})_k$ of velocity and strain-rate fields:
\begin{subequations}\label{eq:sqce-primal}
\begin{itemize}
\item The first idea is to solve the problems in \eqref{eq:sqce-sr} and \eqref{eq:sqce-u} once again, this time with $\bs{\tau}^{(k)}$ instead of $\bs{\hat{\tau}}^{(k)}$:
\begin{equation}\label{eq:sqce-primal-prx}
\begin{dcases}
\bs{\sr}^{(k)} = \argmin_{\bs{\sr} \in Q} \Set{ b(\bs\sr) + j(\bs\sr) - \left\langle\bs{\tau}^{(k)},\bs\sr\right\rangle_{Q*,Q} }\\
\bs{u}^{(k)} = \argmin_{\bs u \in U} \Set{ -\langle \bs f, \bs u\rangle_{\bs L^{r^*}(\Omega),\bs L^r(\Omega)} + \iota_{U_{0\mathrm{D}}}(\bs u) + \frac{1}{2L} \left\lVert \sg \bs u - \left(\bs{\sr}^{(k)} - L \bs{\tau}^{(k)}\right) \right\rVert_Q^2 }
\end{dcases}
\end{equation}
\item The second idea is to obtain $\bs{u}^{(k)}$ from $\bs{\sr}^{(k)}$ by solving the equation $\sg \bs{u}^{(k)} = \bs{\sr}^{(k)}$ in a least-squares sense:
\begin{equation}\label{eq:sqce-primal-lsq}
\begin{dcases}
\bs{\sr}^{(k)} = \argmin_{\bs{\sr} \in Q} \Set{ b(\bs\sr) + j(\bs\sr) - \left\langle\bs{\tau}^{(k)},\bs\sr\right\rangle_{Q*,Q} }\\
\bs{u}^{(k)} = \argmin_{\bs u \in U_{0\mathrm{D}}} \Set{ \frac{1}{2} \left\lVert \sg \bs u - \bs{\sr}^{(k)} \right\rVert_Q^2 }
\end{dcases}
\end{equation}
\item The third idea is to simply set
\begin{equation}\label{eq:sqce-primal-lpt}
\begin{dcases}
\bs{\sr}^{(k)} = \bs{\hat{\sr}}^{(k)}\\
\bs{u}^{(k)} = \bs{\hat{u}}^{(k)}.
\end{dcases}
\end{equation}
\end{itemize}
\end{subequations}

Combining Algorithm \ref{alg:fista}, Lemma \ref{lem:prox} and these definitions of the primal sequence, we obtain the accelerated dual gradient method FISTA* for solving \eqref{eq:min}.

\begin{algorithm}[FISTA* / Accelerated Dual Proximal Gradient Method]\label{alg:dualfista}
\emph{Input:} $\bs\tau^{(0)}\in S$, $\mathrm{gradTol} > 0$\\
\emph{Initialisation:} $k = 1$, $t^{(1)}=1$, $\bs{\hat{\tau}}^{(1)} = \bs\tau^{(0)}$
\begin{description}
\item[(FISTA*.1)\label{it:dualfista1}] Set $L^{(k)} = L$ or find a valid $L^{(k)} > 0$ by backtracking and evaluate
\begin{equation*}
\bs{\hat{\sr}}^{(k)} = \argmin_{\bs \sr \in Q} \Set{ b(\bs\sr) + j(\bs\sr) - \left\langle\bs{\hat{\tau}}^{(k)},\bs\sr\right\rangle_{Q*,Q} }.
\end{equation*}
\item[(FISTA*.2)\label{it:dualfista2}] Solve
\begin{equation*}
\bs{\hat{u}}^{(k)} = \argmin_{\bs u \in U} \Set{ -\langle \bs f, \bs u\rangle_{\bs L^{r^*}(\Omega),\bs L^r(\Omega)} + \iota_{U_{0\mathrm{D}}}(\bs u) + \frac{1}{2L^{(k)}} \left\lVert \sg \bs u - \left(\bs{\hat{\sr}}^{(k)} - L^{(k)} \bs{\hat{\tau}}^{(k)}\right) \right\rVert_Q^2 }.
\end{equation*}
\item[(FISTA*.3)\label{it:dualfista3}] Update
\begin{equation*}
\bs{\tau}^{(k)} = \hat{\bs{\tau}}^{(k)} + \frac{1}{L^{(k)}} \left( \sg \bs {\hat{u}}^{(k)} - \bs{\hat{\sr}}^{(k)} \right).
\end{equation*}
\item[(FISTA*.4)\label{it:dualfista4}] Obtain $\bs u^{(k)}$ and $\bs\sr^{(k)}$ from \eqref{eq:sqce-primal}.
\item[(FISTA*.5)\label{it:dualfista5}] \emph{If} $\left\lVert \sg\bs u^{(k)} - \bs\sr^{(k)} \right\rVert_Q \leq \mathrm{gradTol}$ \emph{Then Return} $\bs u^{(k)}, \bs\sr^{(k)}$ and $\bs\tau^{(k)}$ and \emph{Stop}.
\item[(FISTA*.6)\label{it:dualfista6}] Compute
\begin{equation*}
t^{(k+1)} = \frac{1+\sqrt{1+4t^{(k)2}}}{2}
\end{equation*}
to update the leading point
\begin{equation*}
\bs{\hat{\tau}}^{(k+1)} = \bs\tau^{(k)} + \frac{t^{(k)}-1}{t^{(k+1)}} \left( \bs\tau^{(k)} - \bs\tau^{(k-1)} \right).
\end{equation*}
\item[(FISTA*.7)\label{it:dualfista7}] Set $k \leftarrow k+1$ and \emph{Go To} \ref{it:dualfista1}.
\end{description}
\end{algorithm}

Let us briefly analyse the minimisation problems in \ref{it:dualfista1} and \ref{it:dualfista2}: from properties of convex conjugates it follows immediately \cite[Ch~2]{Treskatis2016} that \ref{it:dualfista1} explicitly reads
\begin{equation*}
\bs{\hat{\sr}}^{(k)} = \nabla F (\bs{\hat{\tau}}^{(k)})
\end{equation*}
with $\nabla F$ from \eqref{eq:dual-grad}. \ref{it:dualfista2} is a variational formulation of the following Stokes problem: find $\bs{\hat{u}}^{(k)}\in U_{0\mathrm{D}}$ and $\hat{p}^{(k)}\in L^{r^*}(\Omega)$ such that
\begin{align*}
-\frac{1}{L^{(k)}}\Div \sg \bs{\hat{u}}^{(k)} + \nabla \hat{p}^{(k)} &= \bs f - \Div \left(\frac{1}{L^{(k)}}\bs{\hat{\sr}}^{(k)} -  \bs{\hat{\tau}}^{(k)}\right)\\
\div \bs{\hat{u}}^{(k)} &= 0.
\end{align*} 
Analogously, the least-squares problem in \eqref{eq:sqce-primal-lsq} leads to the Stokes problem
\begin{align*}
-\Div \sg \bs{u}^{(k)} + \nabla p^{(k)} &= - \Div \bs{\sr}^{(k)}\\
\div \bs{u}^{(k)} &= 0.
\end{align*} 

\subsection{Comparison with Classical Algorithms}

Since conventional algorithms do not include the extrapolation step to determine a leading point, let us consider the non-inertial counterpart of the dual FISTA method, i.e. the dual ISTA method. In the absence of a leading point, defining a separate primal sequence becomes redundant.

\begin{algorithm}[ISTA* / Dual Proximal Gradient Method]\label{alg:dualista}
\emph{Input:} $\bs\tau^{(0)}\in S$, $\mathrm{gradTol} > 0$\\
\emph{Initialisation:} $k = 1$
\begin{description}
\item[(ISTA*.1)\label{it:dualista1}] Set $L^{(k)} = L$ or find a valid $L^{(k)} > 0$ by backtracking and evaluate
\begin{equation*}
\bs{\sr}^{(k)} = \argmin_{\bs \sr \in Q} \Set{ b(\bs\sr) + j(\bs\sr) - \left\langle\bs{\tau}^{(k)},\bs\sr\right\rangle_{Q*,Q} }.
\end{equation*}
\item[(ISTA*.2)\label{it:dualista2}] Solve
\begin{equation*}
\bs{u}^{(k)} = \argmin_{\bs u \in U} \Set{ -\langle \bs f, \bs u\rangle_{\bs L^{r^*}(\Omega),\bs L^r(\Omega)} + \iota_{U_{0\mathrm{D}}}(\bs u) + \frac{1}{2L^{(k)}} \left\lVert \sg \bs u - \left(\bs{\sr}^{(k)} - L^{(k)} \bs{\tau}^{(k)}\right) \right\rVert_Q^2 }.
\end{equation*}
\item[(ISTA*.3)\label{it:dualista3}] Update
\begin{equation*}
\bs{\tau}^{(k)} = \bs{\tau}^{(k-1)} + \frac{1}{L^{(k)}} \left( \sg \bs {u}^{(k)} - \bs{\sr}^{(k)} \right).
\end{equation*}
\item[(ISTA*.4)\label{it:dualista4}] \emph{If} $\left\lVert \sg\bs u^{(k)} - \bs\sr^{(k)} \right\rVert_Q \leq \mathrm{gradTol}$ \emph{Then Return} $\bs u^{(k)}, \bs\sr^{(k)}$ and $\bs\tau^{(k)}$ and \emph{Stop}.
\item[(ISTA*.5)\label{it:dualista5}] Set $k \leftarrow k+1$ and \emph{Go To} \ref{it:dualista1}.
\end{description}
\end{algorithm}

The alternating direction method of multipliers applied to the formulation \eqref{eq:min-split-aug} reads as follows:

\begin{algorithm}[ALG2 / Alternating Direction Method of Multipliers]\label{alg:admm}
\emph{Input:} $\bs\sr^{(0)}\in Q$, $\bs\tau^{(0)}\in S$, a positive sequence $(s^{(k)})_k$ of step sizes, $\mathrm{gradTol} > 0$\\
\emph{Initialisation:} $k = 1$
\begin{description}
\item[(ALG2.1)\label{it:admm1}] Solve
\begin{equation*}
\bs{u}^{(k)} = \argmin_{\bs u \in U} \Set{ -\langle \bs f, \bs u\rangle_{\bs L^{r^*}(\Omega),\bs L^r(\Omega)} + \iota_{U_{0\mathrm{D}}}(\bs u) + \frac{\varrho}{2} \left\lVert \sg \bs u - \left(\bs{\sr}^{(k)} - \frac{1}{\varrho} \bs{\tau}^{(k)}\right) \right\rVert_Q^2 }.
\end{equation*}
\item[(ALG2.2)\label{it:admm2}] Solve
\begin{equation*}
\bs{\sr}^{(k)} = \argmin_{\bs \sr \in Q} \Set{ b(\bs\sr) + j(\bs\sr) - \left\langle\bs{\tau}^{(k)},\bs\sr\right\rangle_{Q*,Q} + \frac{\varrho}{2} \left\lVert \sg \bs u^{(k)} - \bs{\sr} \right\rVert_Q^2 }.
\end{equation*}
\item[(ALG2.3)\label{it:admm3}] Update
\begin{equation*}
\bs{\tau}^{(k)} = \bs{\tau}^{(k-1)} + s^{(k)} \left( \sg \bs {u}^{(k)} - \bs{\sr}^{(k)} \right).
\end{equation*}
\item[(ALG2.4)\label{it:admm4}] \emph{If} $\left\lVert \sg\bs u^{(k)} - \bs\sr^{(k)} \right\rVert_Q \leq \mathrm{gradTol}$ \emph{Then Return} $\bs u^{(k)}, \bs\sr^{(k)}$ and $\bs\tau^{(k)}$ and \emph{Stop}.
\item[(ALG2.5)\label{it:admm5}] Set $k \leftarrow k+1$ and \emph{Go To} \ref{it:admm1}.
\end{description}
\end{algorithm}

We observe a number of similarities and differences:
\begin{itemize}
\item The steps \ref{it:admm1} and \ref{it:dualista2} are identical, provided that $\varrho = 1/L^{(k)}$.
\item The steps \ref{it:admm2} and \ref{it:dualista1} are identical, provided that $\varrho = 0$. It can be seen as a major disadvantage of the alternating direction method of multipliers with $\varrho > 0$ that for Casson and Herschel-Bulkley fluids, the solution of \ref{it:admm2} cannot be obtained explicitly. Instead, an iterative scheme such as Newton's method is required in every iteration of ALG2 to compute $\bs{\sr}^{(k)}$. In \ref{it:dualista1}, ISTA* only evaluates a nonlinear function. The same holds true for \ref{it:dualfista1}.
\item The updates \ref{it:admm3}-\ref{it:admm5} and \ref{it:dualista3}-\ref{it:dualista5} are identical, provided that $s^{(k)} = 1/L^{(k)}$.
\item The parameters $\varrho$ and $s^{(k)}$ in ALG2 have to be chosen heuristically, whereas a globally optimal value for $L^{(k)}$ can be calculated analytically, or estimated by backtracking.
\end{itemize}

In conclusion, ISTA* solves either the same or simpler subproblems than ALG2 in every iteration. FISTA* additionally incorporates an extrapolation step to achieve a higher rate of convergence. We will quantify these rates for each algorithm in the next subsection.

\subsection{Convergence of FISTA*}

In step \ref{it:dualfista4}, \eqref{eq:sqce-primal} provides three alternatives for defining the iterates $\bs u^{(k)}$ and $\bs \sr^{(k)}$. FISTA* with \eqref{eq:sqce-primal-prx} or \eqref{eq:sqce-primal-lsq} approximately doubles the computational cost per iteration compared to ISTA*: the lion share of computing power is required for the execution of \ref{it:dualfista2} and, to a far lesser extent, \ref{it:dualfista1}. The additional step \ref{it:dualfista4} demands for a solution of these two problems with different data a second time. In contrast, if the primal sequence \eqref{eq:sqce-primal-lpt} based on the leading point $\bs{\hat{\tau}}^{(k)}$ is chosen, then \ref{it:dualfista4} becomes trivial and no extra problems have to be solved.

Unfortunately, convergence of the sequence \eqref{eq:sqce-primal-lpt} cannot be guaranteed. At this stage, it is an open problem whether the convergence of this sequence might actually fail, whether it generally converges even though no proof appears to be available, or if additional assumptions on the problem are required to ensure convergence (cf \cite[Ch~3]{Treskatis2016}).

For the other two alternatives, we have the following result:

\begin{theorem}[Convergence of FISTA* and ISTA*]
Let $\bs{\bar{u}}$ be the exact solution of Problem \eqref{eq:min}.

If the sequences $(\bs u^{(k)})_k$, $(\bs \sr^{(k)})_k$ and $(\bs \tau^{(k)})_k$ are generated by FISTA* with either \eqref{eq:sqce-primal-prx} or \eqref{eq:sqce-primal-lsq}, then
\begin{align*}
\left\lVert \bs\sr^{(k)} - \sg \bs{\bar{u}} \right\rVert_Q &= O(1/k)\\
\left\lVert \bs u^{(k)} - \bs{\bar{u}} \right\rVert_U &= O(1/k).
\end{align*}

If the sequences $(\bs u^{(k)})_k$, $(\bs \sr^{(k)})_k$ and $(\bs \tau^{(k)})_k$ are generated by ISTA*, then
\begin{align*}
\left\lVert \bs\sr^{(k)} - \sg \bs{\bar{u}} \right\rVert_Q &= O(1/\sqrt{k})\\
\left\lVert \bs u^{(k)} - \bs{\bar{u}} \right\rVert_U &= O(1/\sqrt{k}).
\end{align*}
\begin{proof}
We refer to \cite[Ch~3]{Treskatis2016}.
\end{proof}
\end{theorem}

Let us mention once again that under additional regularity assumptions on the problem, the alternating direction method of multipliers ALG2 with a suitable choice of parameters exhibits the same convergence rate as ISTA.

The sequence of dual multipliers $(\bs \tau^{(k)})_k$  is bounded and therefore possesses a weakly convergent subsequence. If weak convergence of the entire sequence $(\bs \tau^{(k)})_k$ is desired, then this can be achieved by a modification of the extrapolation sequence $(t^{(k)})_k$ in FISTA* \cite{Chambolle2015}.

\subsection{Computational Techniques}

\paragraph{Solution of the Stokes Problems}

For solving the Stokes problems that occur in each method, we apply the classical preconditioned conjugate gradient Uzawa (PCGU) method of Cahouet and Chabard \cite[pp~892--893]{Cahouet1988}. Glowinski \cite[Sec~20--22]{Glowinski2003} motivates this method by successively improving on the very basic Uzawa (i.e. dual gradient) method for the Stokes problem. These step-by-step improvements on the speed of convergence of the algorithm are achieved by
\begin{itemize}
\item using conjugate gradients instead of steepest descent.
\item performing an exact line search at every iteration instead of using the global Lipschitz constant of the dual gradient as a worst-case estimate. Since the problem is quadratic, the exact step size is straightforward to calculate.
\item preconditioning the problem. While this is essential for the instationary generalisation of the Stokes problem, it is not beneficial for the stationary case that we consider here. With no time dependence, the suggested preconditioner would simply degenerate to a multiple of the identity. Looking at the discrete problem, one could interpret the occurrence of a mass matrix as preconditioning, though.
\end{itemize}

We terminate the PCGU algorithm as soon as the dual gradient is sufficiently small, as measured by a positive constant $\mathrm{stokesTol}$:
\begin{equation*}
\left\lVert \div \bs u^{(k)} \right\rVert_{L^r(\Omega)} \leq \mathrm{stokesTol}.
\end{equation*}

Our convergence analysis of the dual proximal gradient methods assumes that the proximal map is evaluated exactly. Therefore, we should at least set $\mathrm{stokesTol} \ll \mathrm{gradTol}$. This poses no major obstacle, since the PCGU algorithm achieves a linear rate of convergence for the Stokes problem, compared to the sublinear rate of the outer optimisation loop.

\paragraph{Adaptive Re-Starting}

In general, FISTA is not a monotone method in the sense that the value of the dual objective $F+G$ may also increase from one iteration to another. Analogously, the primal sequence generated by Algorithm FISTA* may temporarily digress from the solution $(\bs{\bar{u}},\bs{\bar{\sr}})$. This is a well-known property of accelerated gradient schemes and can be interpreted as excessive momentum from past iterations, that causes the sequence to overshoot the minimiser and to converge in a spiralling motion.

In \cite{Beck2009a}, Beck and Teboulle propose a simple modification of FISTA, called MFISTA, which guarantees monotonicity. In contrast to the basic FISTA scheme, it requires the functional $F+G$ to be evaluated at every iteration.

O'Donoghue and Candes \cite{ODonoghue2013} suggest to adaptively re-start the algorithm once an increase in the objective is detected in order to preserve monotonicity and to discard accumulated momentum of the iteration. Rather than observing the functional values, they showed a re-starting criterion based solely on the dual gradient to be similarly effective. According to this \emph{gradient scheme}, Algorithm FISTA* is re-started whenever
\begin{equation}\label{eq:restart}
\left\langle\sg \bs{\hat{u}}^{(k)} - \bs{\hat{\sr}}^{(k)},\bs \tau^{(k)}-\bs\tau^{(k-1)}\right\rangle_{S*,S} < 0.
\end{equation}
In that case, $\bs \tau^{(k)}-\bs\tau^{(k-1)}$ would be an ascent direction for the dual functional at $\bs{\hat{\tau}}^{(k)}$. The authors of \cite{ODonoghue2013} point out that this scheme unites the benefits of increased numerical stability near the optimum on the one hand and, on the other hand, no extra computational expenditure: all quantities in \eqref{eq:restart} have already been computed previously.

By allowing for re-starts, the worst-case convergence rate decreases from $O(1/k)$ to $O(1/\sqrt{k})$. This is a consequence of the fact that the first step after (re-)starting is equivalent to a step in the unaccelerated dual proximal gradient method ISTA* and we have
\begin{equation*}
\bs\tau^{(k)} = \bs\tau^{(k-1)} + \frac{1}{L^{(k)}} (\sg \bs{\hat{u}}^{(k)} - \bs{\hat{\sr}}^{(k)}).
\end{equation*}
Therefore,
\begin{equation*}
\left\langle\sg \bs{\hat{u}}^{(k)} - \bs{\hat{\sr}}^{(k)},\bs \tau^{(k)}-\bs\tau^{(k-1)}\right\rangle_{S*,S} = \frac{1}{L^{(k)}} \left\lVert \sg \bs{\hat{u}}^{(k)} - \bs{\hat{\sr}}^{(k)} \right\rVert_Q^2 \nless 0.
\end{equation*}
and this first step is generally accepted.

Surely, the rationale behind re-starting schemes is that re-initialisations only occur as isolated events, not after every single iteration. This way, the convergence rate would remain close to $O(1/k)$. A re-start would ideally result in a shortcut towards the solution and thus decrease the error more efficiently than continued iterations with full momentum. For a numerical study of this effect, we refer to \cite{ODonoghue2013}.

\paragraph{Discretisation}

We assume that $\Omega$ is a polygonal domain we let $\mathcal{T}_h$ be a regular triangulation on $\Omega$. We construct a finer triangulation $\mathcal{T}_{h/2}$ by connecting the edge midpoints in each triangle $T\in \mathcal{T}_h$.

Following Glowinski \cite[p~303]{Glowinski2003}, we apply the $\mathds{P}_1$-iso-$\mathds{P}_2/\mathds{P}_1$ element of Bercovier and Pironneau \cite{Bercovier1979} for discretising the velocity and pressure, respectively. $\mathds{P}_k$ denotes the space of polynomials in two variables of degree at most $k$. We approximate the strain rate and stress with piecewise constant elements on the finer mesh.

Overall, we consider the following finite-element spaces:
\begin{align*}
U_{*0,h} &:= \Set{ \bs u_h \in C(\bar{\Omega})^2 | \bs u_h\rvert_\Gamma = 0 \quad \text{and} \quad \bs u_h\rvert_T \in \mathds{P}_1^2, \quad \forall T \in \mathcal{T}_{h/2}  }\\
P_{h} &:= \Set{ p_h \in C(\bar{\Omega}) | p_h \rvert_T \in \mathds{P}_1, \quad \forall T \in \mathcal{T}_h }\\
Q_h &:= \Set{ \bs\sr_h \in L^r(\Omega)^3 | \bs\sr_h \rvert_T \in \mathds{P}_0^3, \quad  \forall T \in \mathcal{T}_{h/2} }\\
S_h &:= \Set{ \bs\tau_h \in L^{r^*}(\Omega)^3 | \bs\tau_h \rvert_T \in \mathds{P}_0^3, \quad  \forall T \in \mathcal{T}_{h/2} }.
\intertext{With $\bs u_{\mathrm{D},h}\in C(\Gamma)\cap U_\mathrm{D}$ serving as an approximation to $u_\mathrm{D}$ that is linear on the triangle edges of $\mathcal{T}_{h/2}$,we also introduce the convex set}
U_{*\mathrm{D},h} &:= \Set{ \bs u_h \in C(\bar{\Omega})^2 | \bs u_h\rvert_\Gamma = \bs u_{\mathrm{D},h} \quad \text{and} \quad \bs u_h\rvert_T \in \mathds{P}_1^2, \quad \forall T \in \mathcal{T}_{h/2}  }.
\end{align*}

It is well-known that these spaces are LBB-stable \cite{Bercovier1979}, i.e. the Stokes problems discretised with these finite elements possess a unique solution.

\section{Numerical Results}
\label{sec:results}

In this section, we will conduct some numerical experiments. Our objective is twofold: firstly, we investigate the effect of the acceleration to demonstrate the improved rate of convergence of the accelerated dual proximal gradient method FISTA* compared to ISTA* and in particular the alternating direction method of multipliers ALG2, the current benchmark for solving viscoplastic flow problems with no regularisation. Secondly, we wish to verify that FISTA* computes accurate approximations, in particular with respect to predicting yielded and unyielded regions in the flow. A simple scheme for adaptive mesh refinements will help us to obtain high-fidelity approximations with great efficiency.

We leave a few parameters unaltered for our simulations: for the dual gradient methods, we set $L^{(0)} = 1/2$. This Lipschitz parameter is kept constant for Bingham and Casson flow simulations, while we choose a magnifying factor of $\eta = 1.1$ for the backtracking procedure in the Herschel-Bulkley setting. For ALG2, we use the corresponding value $\varrho = s^{(k)} = 2$ for all $k$. We initialise the PCGU algorithm for solving each Stokes problem with the converged solution for the pressure of the previous run. All other initial guesses for the required variables shall be $0$ in each example. In the stopping criteria we set $\mathrm{gradTol} = 10^{-6}$ and $\mathrm{stokesTol} = 10^{-12}$. We run our programs in \textsf{MATLAB} R2013a 64-bit on a laptop with \textsf{Intel\textsuperscript\textregistered Core\texttrademark i7} CPU 4x2.50 GHz.

For our numerical experiments, we consider two different flow problems in a square reservoir: a force-driven and a wall-driven flow. In both cases, we define $\Omega := \left] 0,1 \right[^2$. To mesh the geometry, we proceed as follows: first, we generate a uniform grid of $1/h \times 1/h$ squares. Then we divide each square diagonally into four congruent triangles, the collection of which defines the coarse pressure grid $\mathcal{T}_h$.

In all our experiments with FISTA*, the primal sequence \eqref{eq:sqce-primal-lpt} based on the leading point turns out to converge. Hence, there is no need to fall back to any of the computationally expensive alternatives \eqref{eq:sqce-primal-prx} or \eqref{eq:sqce-primal-lsq}.

\subsection{Force-Driven Cavity}

In \cite{Reyes2012}, De los Reyes and Gonz{\'a}lez Andrade simulate the flow of a Bingham fluid, which is driven by the force
\begin{equation*}
\bs f (x_1,x_2) := 300 (x_2-0.5,0.5-x_1)^\top
\end{equation*}
with Bingham number $\bi = 10\sqrt{2}$ and homogeneous boundary conditions. We use the same parameters, and additionally consider the corresponding Casson and Herschel-Bulkley flow problems ($r = 1.5$). Our simulations are carried out on the grid with $h = 1/32$.

\paragraph{Convergence Rates}

We carry out 5,000 iterations of FISTA* to approximate the exact solution of \eqref{eq:min} $\bs{\bar{u}} \approx \bs u_h^{(5,000)}$. This allows us to report very accurate estimates of the true error $\lVert \bs u_h^{(k)} - \bs{\bar{u}}\rVert_U$.

\begin{figure}
	\includegraphics[width=0.49\textwidth]{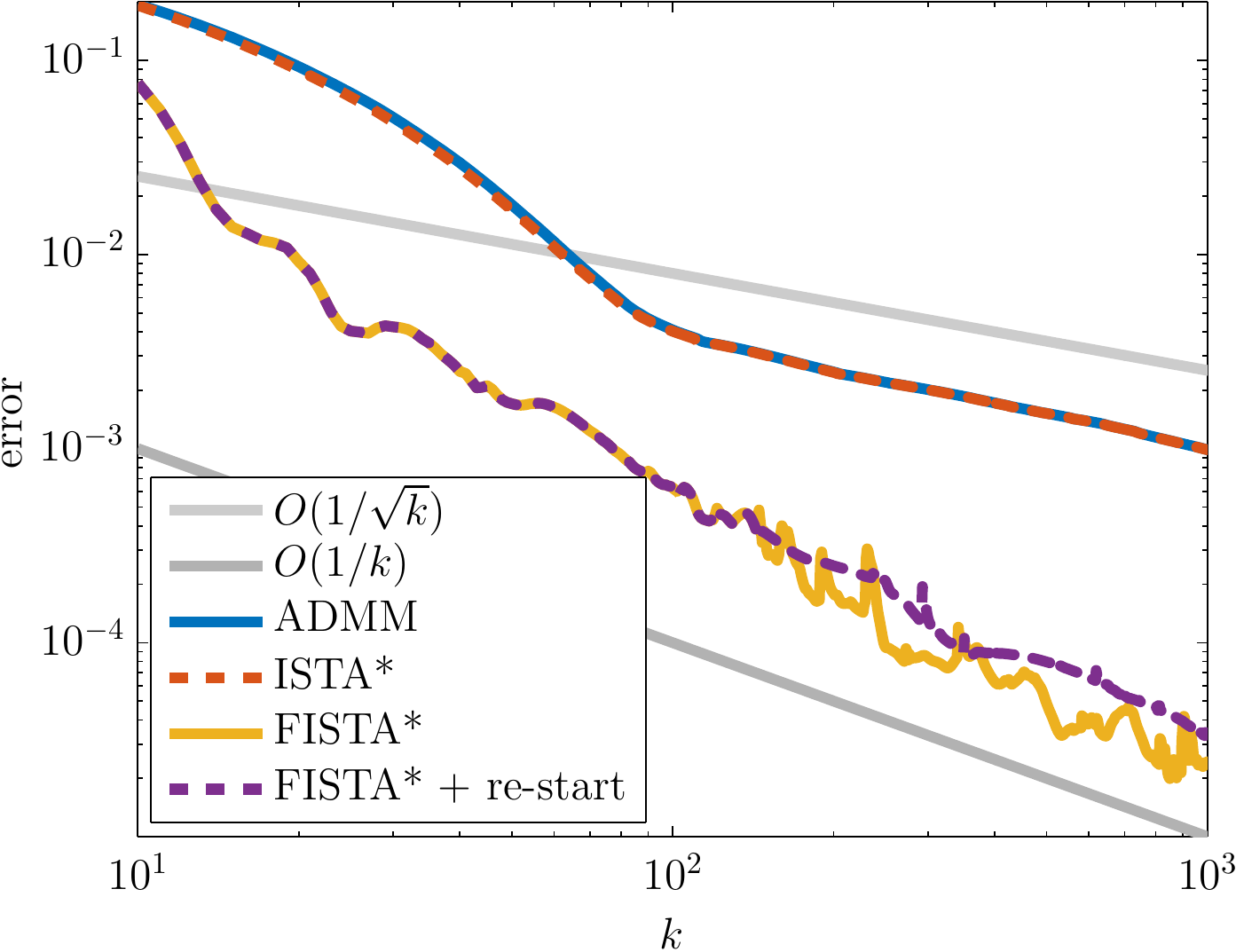}
	\hfill
	\includegraphics[width=0.49\textwidth]{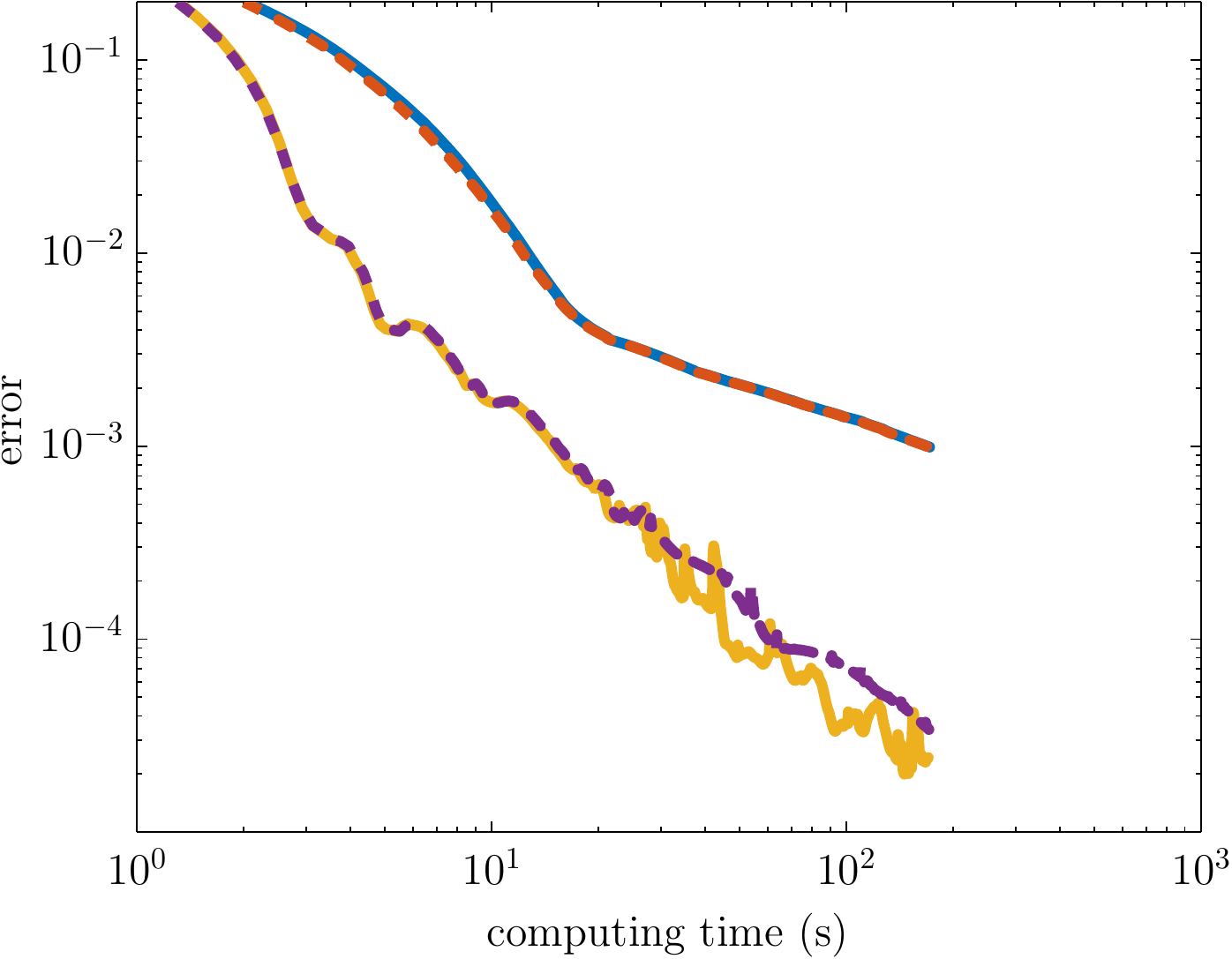}
	\caption{Convergence history for Bingham flow in a force-driven square reservoir.}\label{fig:cvgce-fdc-bi}
\end{figure}

In Figure \ref{fig:cvgce-fdc-bi}, we compare the convergence of ADMM / ALG2, ISTA* and FISTA* for the Bingham flow problem. Additionally, we show how the convergence of FISTA* is affected when the criterion \eqref{eq:restart} is monitored to trigger adaptive re-starts of the method.

\begin{figure}
	\includegraphics[width=0.49\textwidth]{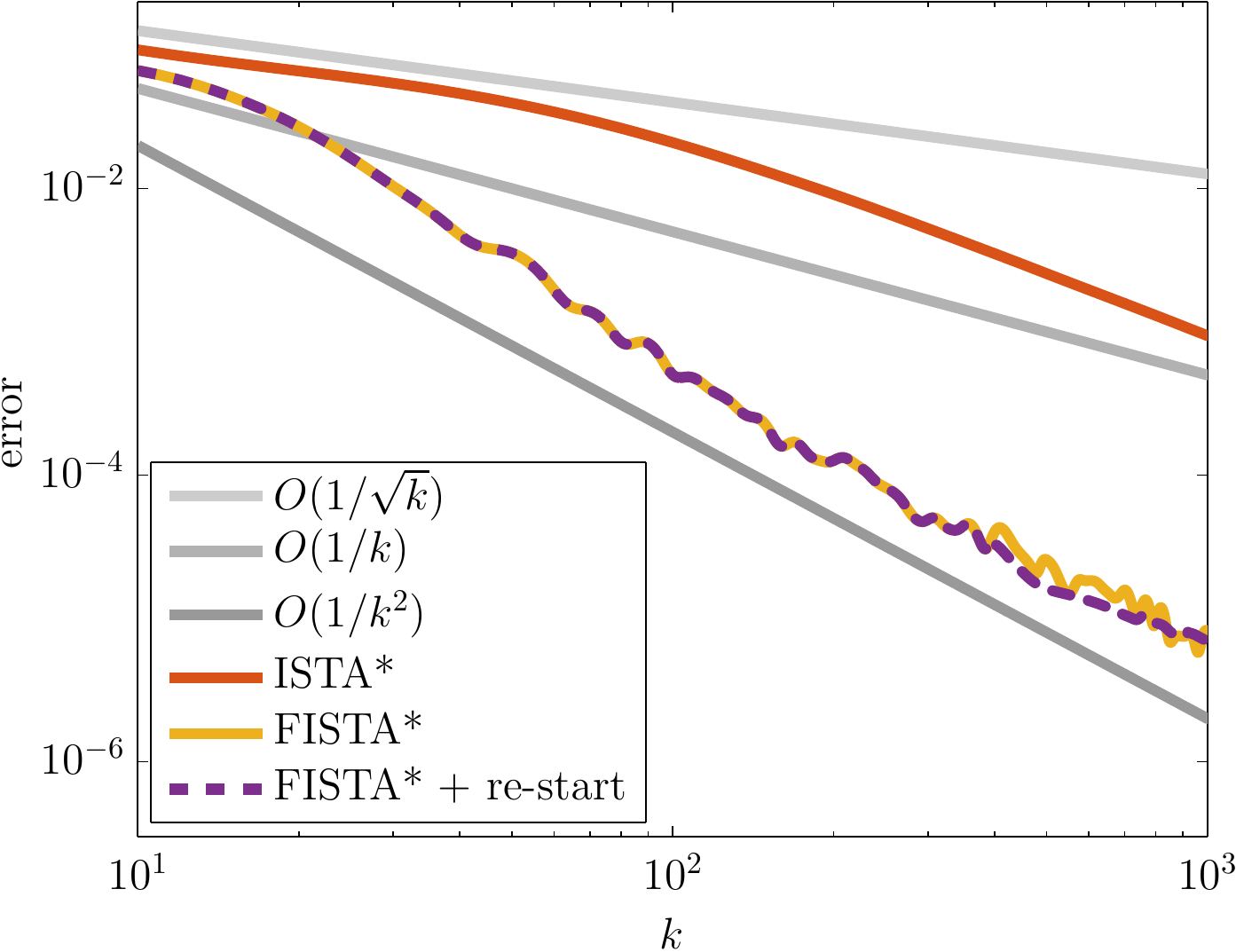}
	\hfill
	\includegraphics[width=0.49\textwidth]{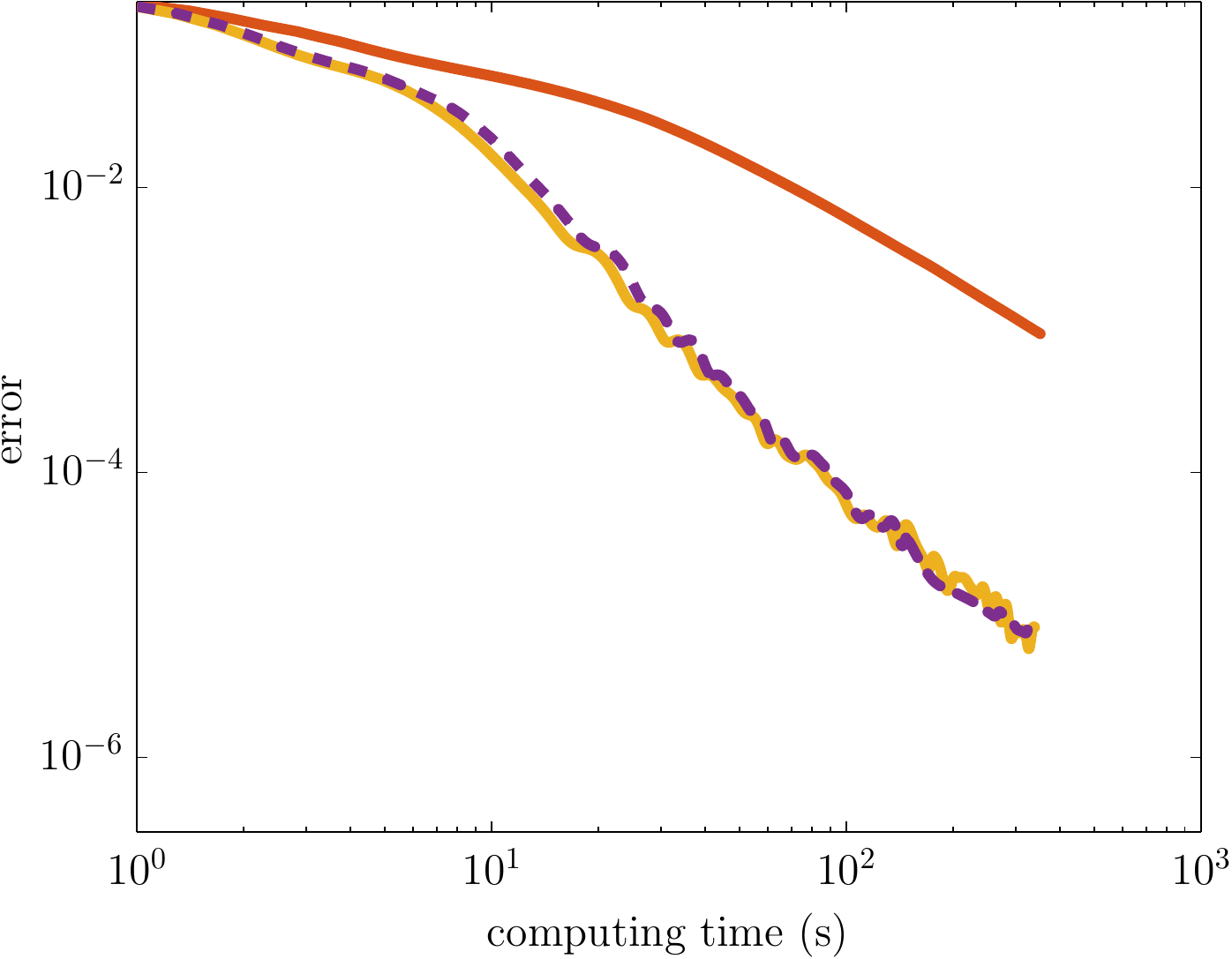}
	\caption{Convergence history for Casson flow in a force-driven square reservoir.}\label{fig:cvgce-fdc-ca}
\end{figure}

\begin{figure}
	\includegraphics[width=0.49\textwidth]{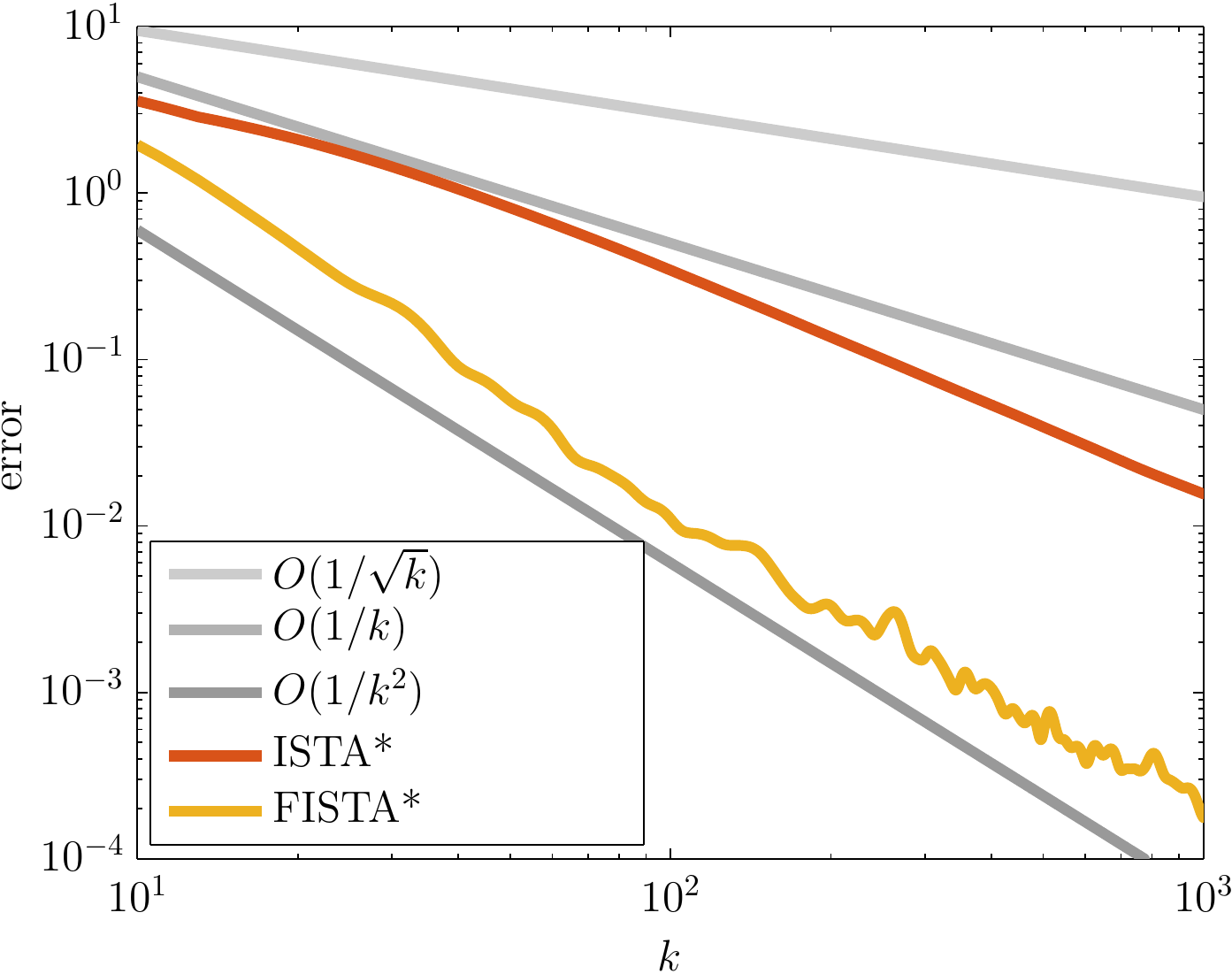}
	\hfill
	\includegraphics[width=0.49\textwidth]{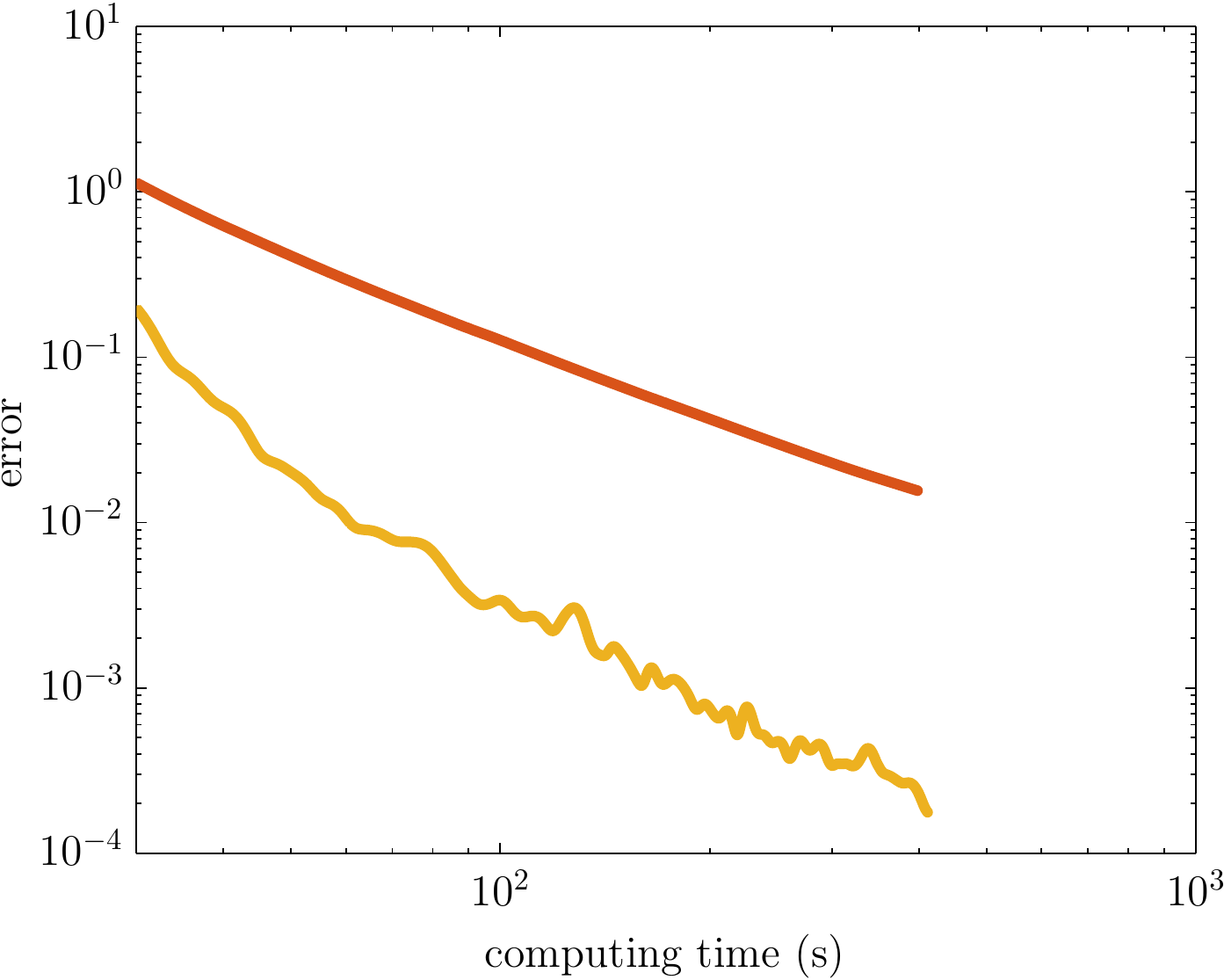}
	\caption{Convergence history for Herschel-Bulkley flow in a force-driven square reservoir ($r = 1.5$).}\label{fig:cvgce-fdc-hb}
\end{figure}

In this setting, the accelerated methods achieve a convergence rate of order $O(1/k)$ within only few iterations. ADMM and ISTA* also exhibit their worst-case convergence rate of $O(1/\sqrt{k})$ after a start-up phase. FISTA* with adaptive re-starting discards its previous momentum after the iterations $k = 144$ and $k = 351$. Although the descent is qualitatively more consistent, the method is not strictly monotonous and in absolute numbers, the errors lie above those that can be achieved without re-starting.

It is a very important observation that replacing iteration numbers with computing time on the horizontal axis does not visually affect the graphs of the error curves. This confirms the theoretical expectations that the cost of one iteration in ADMM, ISTA* or FISTA* is virtually identical. However, after 1,000 iterations, FISTA* has computed an approximation that is about two orders of magnitudes more accurate than the estimates returned by ADMM and ISTA*. Hence, even though FISTA* converges in a non-monotone fashion, it is hard to think of any practical disadvantages of this effect.

Moving on to the results for the Casson fluid in Figure \ref{fig:cvgce-fdc-ca}, fast optimisation algorithms show convergence of order $O(1/k^2)$. Surprisingly, even ISTA* appears to attain the same rate asymptotically. In contrast, the results suggest that re-starts of FISTA* after the iterations with indices $k = 313$ and $k = 919$ decelerate the convergence rate down to $O(1/k)$. Nevertheless, the re-starting scheme effectively improves the monotonicity of the descent. After 1,000 iterations, FISTA* has arrived at a solution about ten times more accurate than the 1,000\textsuperscript{th} iterate of ISTA*. ADMM is uneconomical for solving the Casson flow problem as it requires an accurate solution of a nonlinear problem at every iteration.

In Figure \ref{fig:cvgce-fdc-hb}, we present the results for the simulations of a shear-thinning Herschel-Bulkley fluid. No re-starts of FISTA* occurred during the first 1,000 iterations. The results are otherwise similar to the Casson flow problem.

\paragraph{Yielded and Unyielded Flow Regions}

Let us now compare the geometry of the stagnant zones as they are predicted by the genuinely nonsmooth methods of this article and the regularised approach of De los Reyes and González Andrade in \cite{Reyes2012}. Even though the authors study a time-dependent Bingham flow problem, they report that a steady state is quickly attained. This allows us to compare their results under the quasi-stationary regime with ours.

To visualise the regions of yielded and unyielded flow, it is common to plot areas where $\lvert \bs\tau_h \rvert \leq \bi$ and $\lvert \bs\tau_h \vert > \bi$ in two different colours. However, in its basic form, this approach does not normally provide satisfactory results, as numerical errors near the interface become visible in the form of colourful noise. A traditional remedy \cite{Reyes2010,Yu2007,Olshanskii2009} exploits the fact that the solution tends to be more regular and converges a lot faster in yielded flow regions and therefore the shape of sets with $\lvert \bs\tau_h \rvert \leq (1+\varepsilon)\bi$ for a positive „correction factor“ $\varepsilon$  is typically far smoother. Nevertheless, the main drawback of this postprocessing step is that it introduces a systematic error by overestimating the actual unyielded regions. The extra effort we have invested in the numerical solution by means of a genuinely nonsmooth method would partially be futile if we still apply some form of smoothing to the results.

\begin{figure}
	\includegraphics[width=0.45\textwidth]{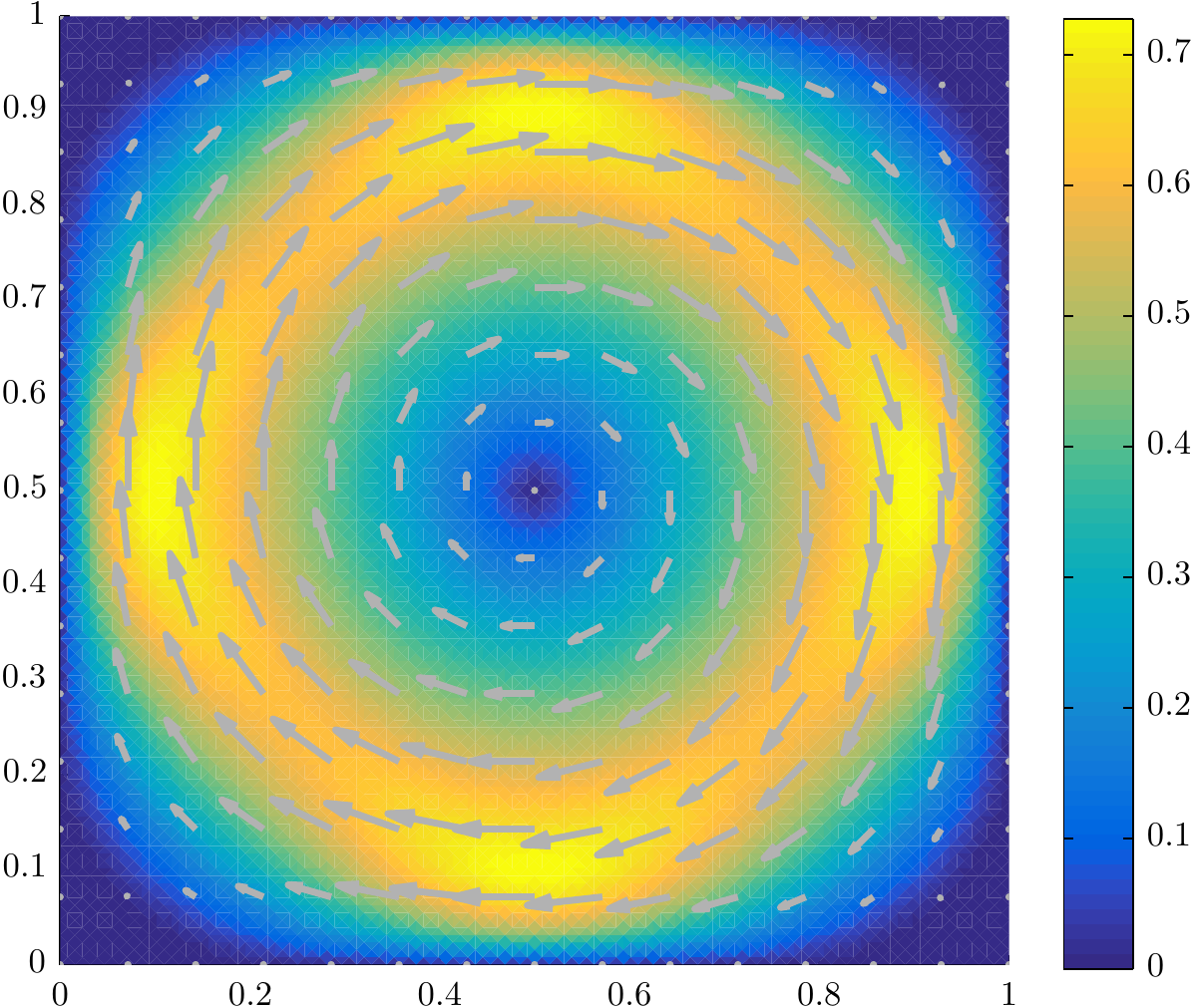}
	\hfill
	\includegraphics[width=0.45\textwidth]{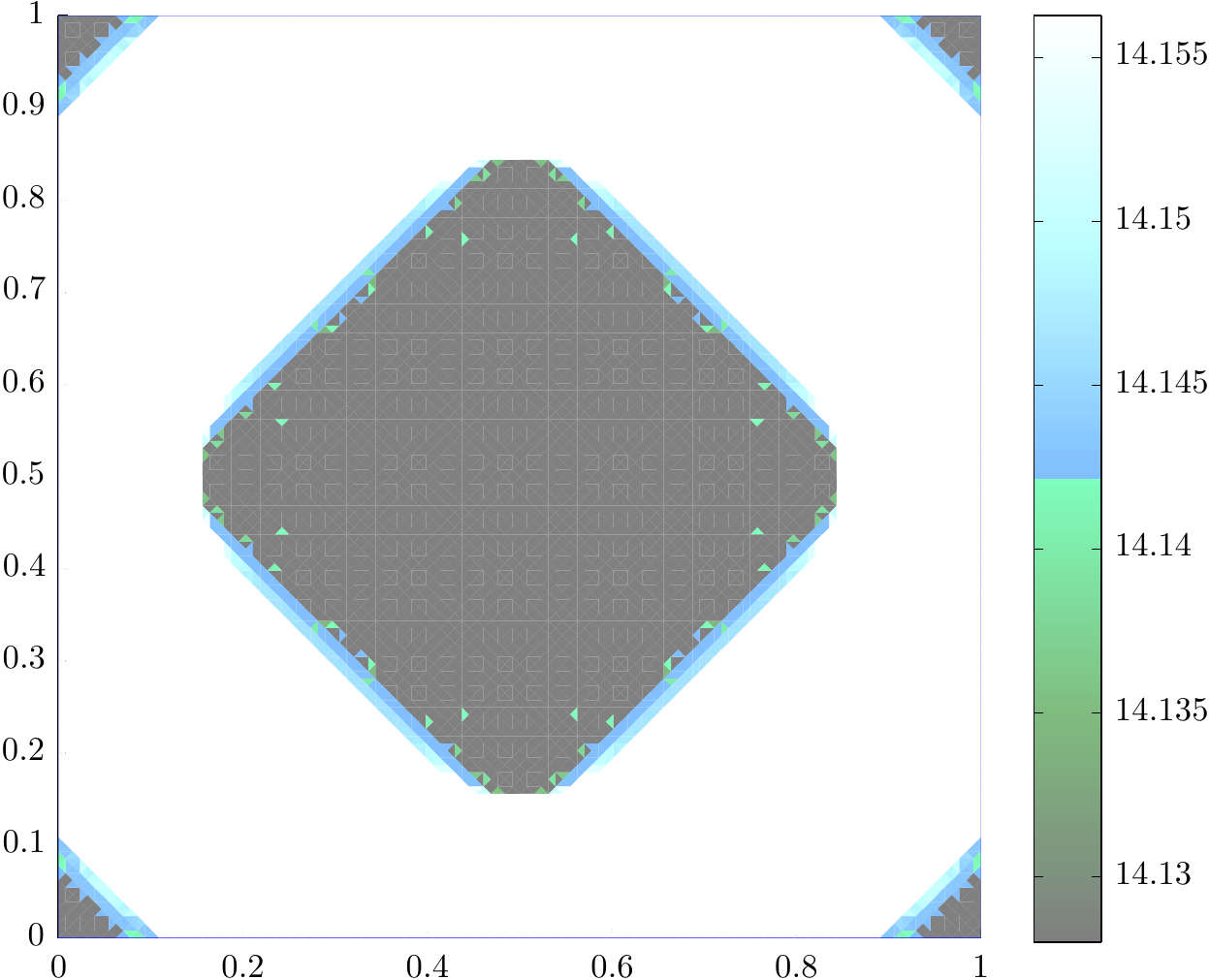}
	\caption{Flow velocity and plug zones for rotational Bingham flow in a square reservoir ($h=1/32$).}\label{fig:fdc-bi}
\end{figure}

\begin{figure}
	\includegraphics[width=0.45\textwidth]{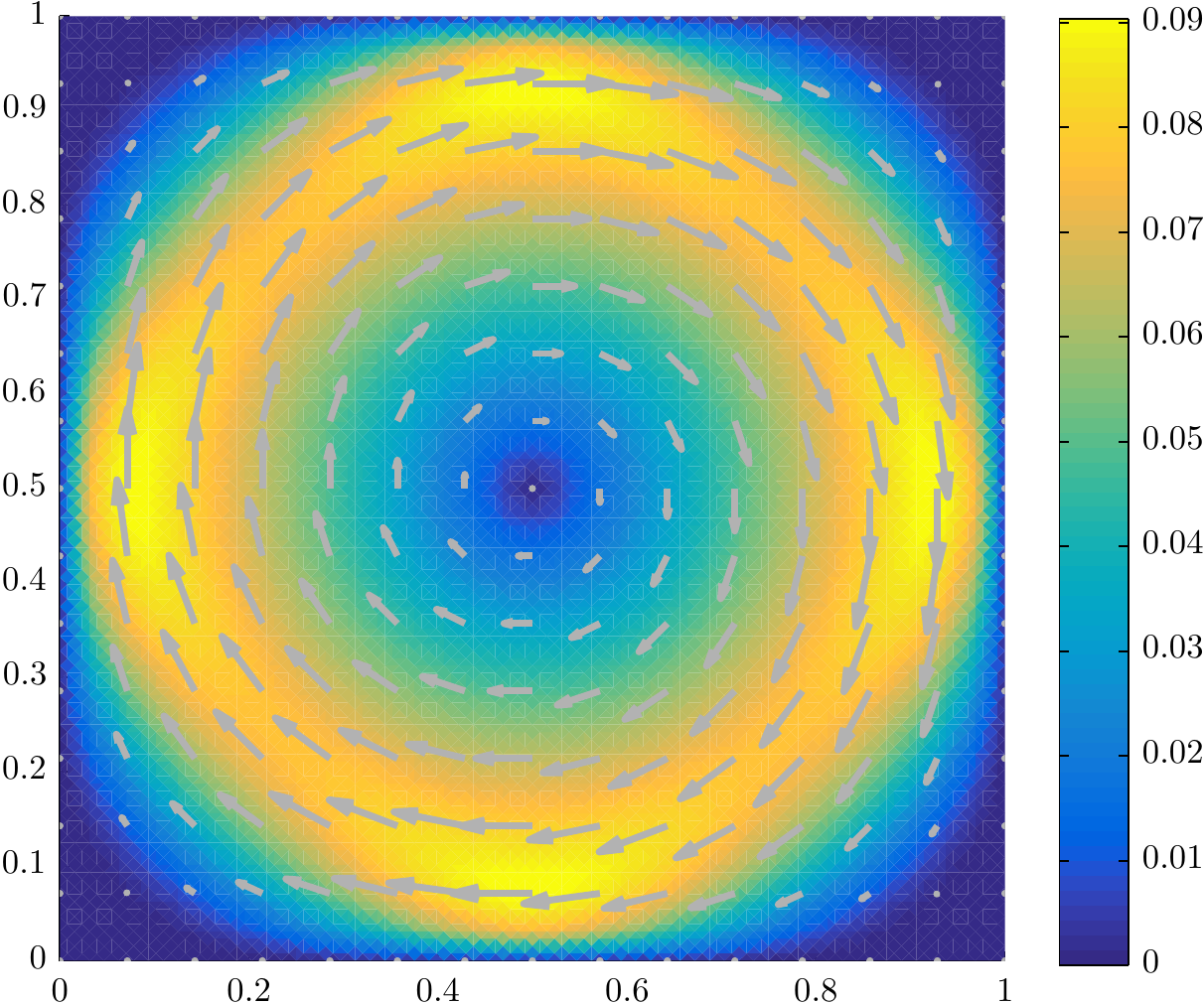}
	\hfill
	\includegraphics[width=0.45\textwidth]{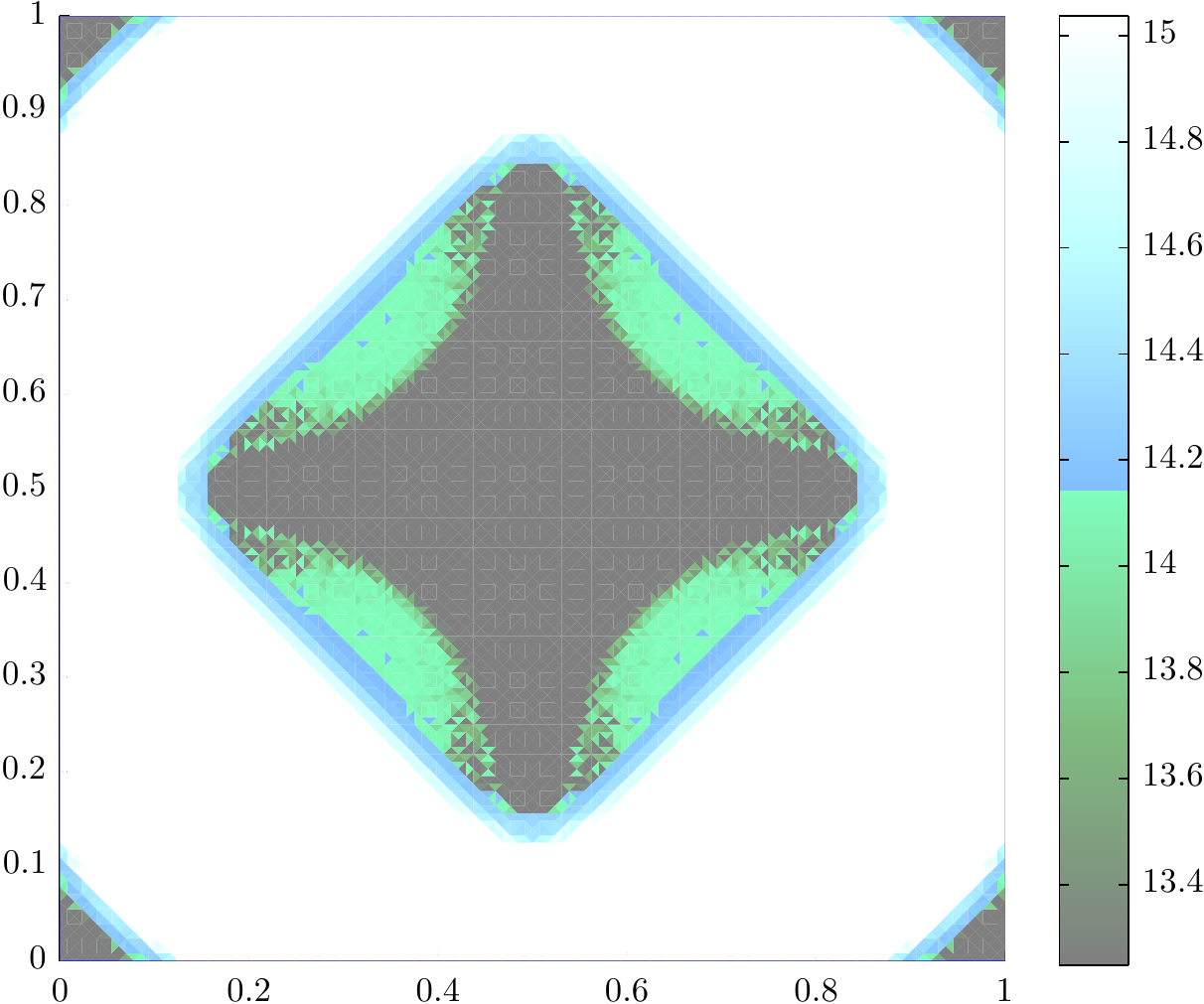}
	\caption{Flow velocity and plug zones for Casson flow.}\label{fig:fdc-ca}
\end{figure}

\begin{figure}
	\includegraphics[width=0.45\textwidth]{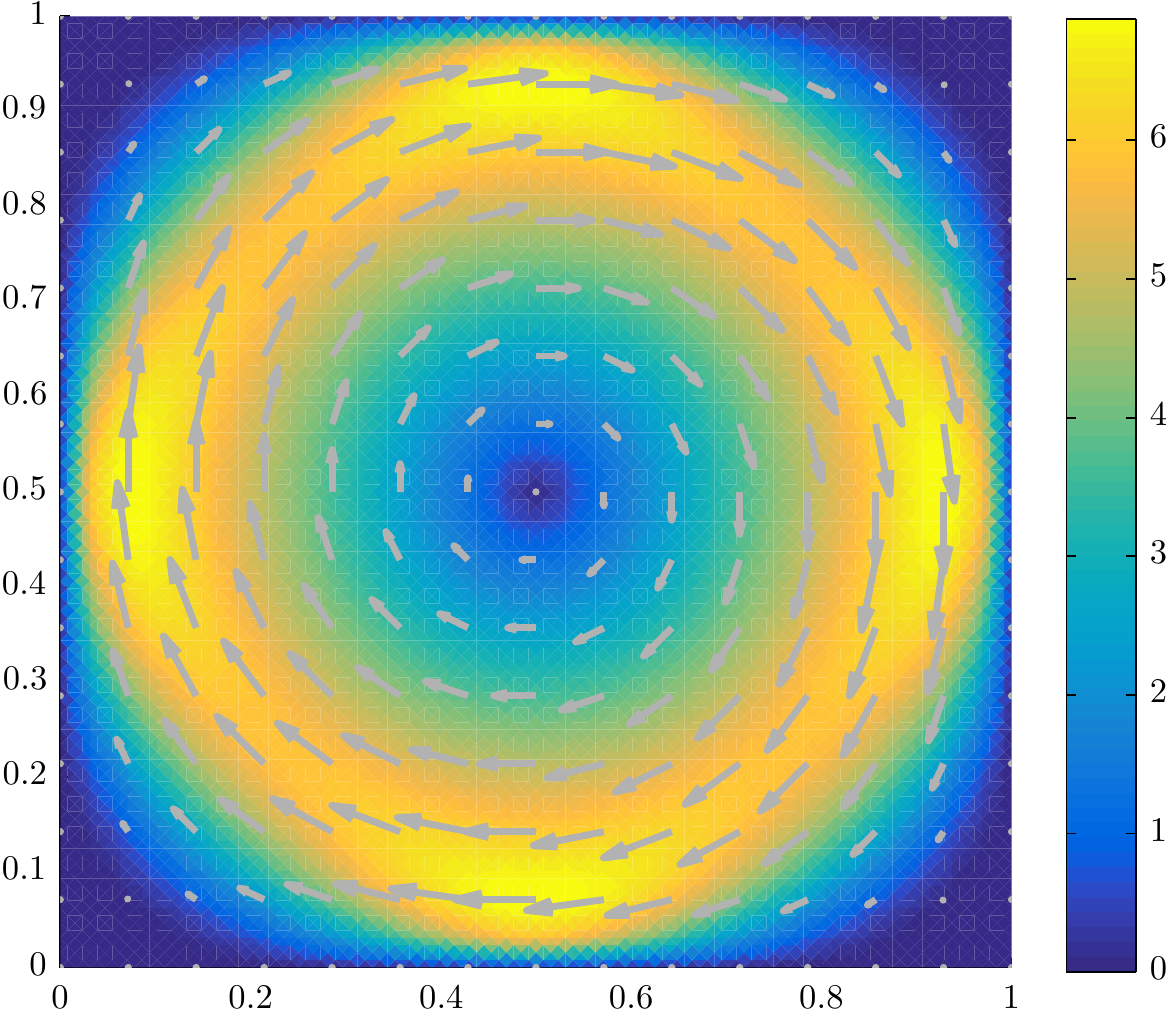}
	\hfill
	\includegraphics[width=0.45\textwidth]{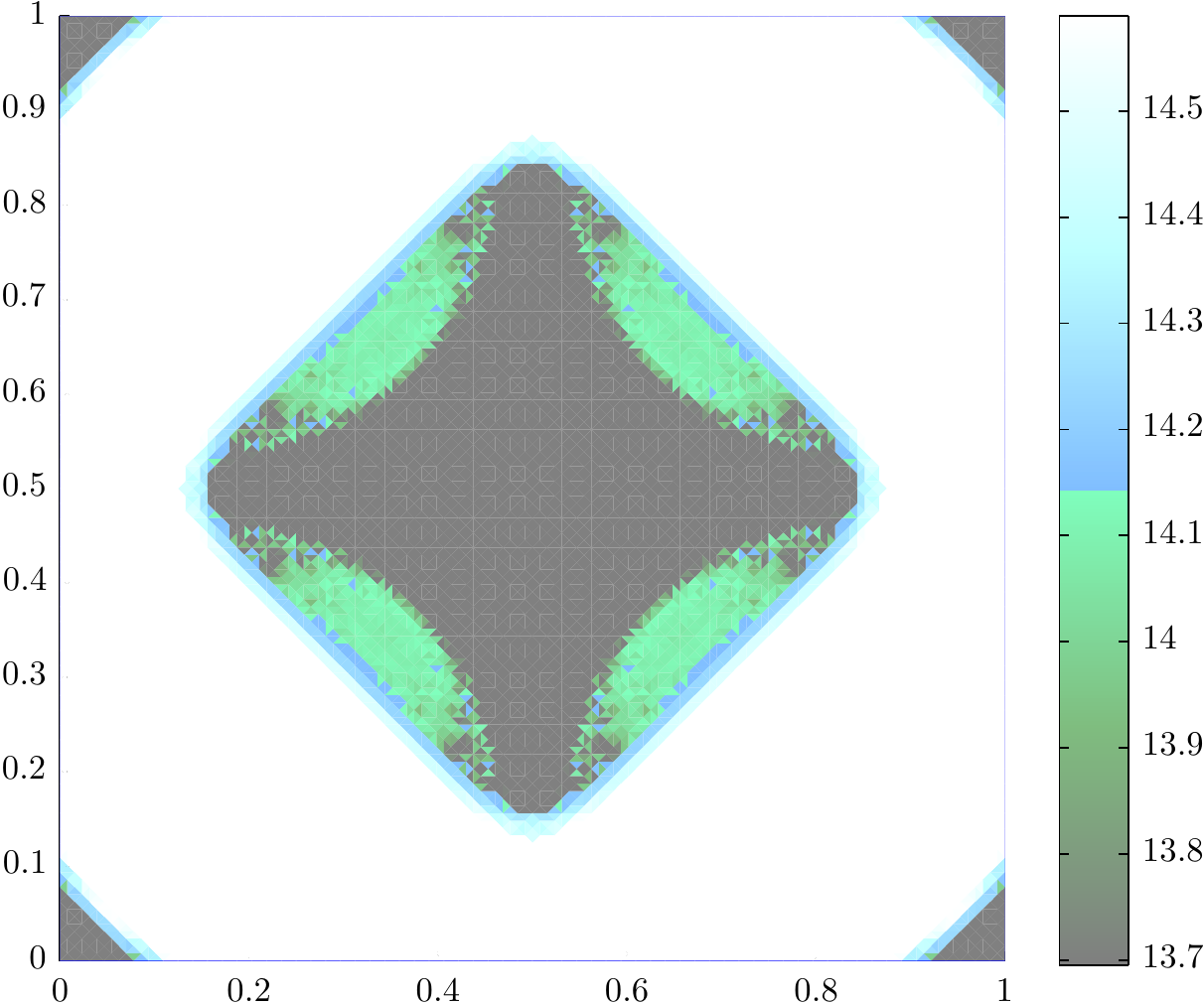}
	\caption{Flow velocity and plug zones for Herschel-Bulkley flow ($r=1.5$).}\label{fig:fdc-hb}
\end{figure}

We therefore propose another means of visualising the yield surfaces: in the right-hand diagrams of Figures \ref{fig:fdc-bi}-\ref{fig:fdc-hb}, we plot the magnitude of the extra stress tensor $\lvert \bs\tau_h \rvert$ in a window of $\pm 0.1 \%$ around the Bingham number $\bi$.  Stress magnitudes below the critical value appear in grey-green, yielded regions are displayed as blue-white. Therefore, the interface between yielded and unyielded regions, as predicted by the numerical solution, lies at the sharp transition from blue to green. Meanwhile, since the classification into „yielded“ and „unyielded“ is least reliable near $\lvert \bs\tau_h \rvert = \bi$ due to numerical errors, the width of blue and green shaded areas serves as an indicator of uncertainty in the identification of flow regions. As for the „correction factor“ $\varepsilon$, there is of course some arbitrariness in choosing the width of the interval around the yield stress, which defines the span of the colour gradients. However, this width never introduces any systematic errors into the visualisation, as the discontinuity in our colourbar always occurs exactly at the value $\bi$.


Upon comparing the results in Figure \ref{fig:fdc-bi} with \cite{Reyes2012}, we observe major differences. The central solid region of approximately square shape deviates significantly from the cross-like or „+“-shape computed by De los Reyes and Gonz{\'a}lez Andrade. In our visualisation, it appears that the stress magnitude lies well below the yield stress, which indicates that those results should be reliable.

In the corresponding graphs for Casson flow (Figure \ref{fig:fdc-ca}) and Herschel-Bulkely flow (Figure \ref{fig:fdc-hb}) we recognise a similar pattern of a rounded cross-like structure in the stress. Still, our numerical results suggest that the actual unyielded region extends beyond the grey cross and in fact the flow stagnates in the entire region that is shaped like a rounded square.

On another hand, the flow field in the Bingham case computed by FISTA* agrees both qualitatively and quantitatively with the semismooth Newton method in \cite{Reyes2010}. This observation confirms that regularisation techniques are most appropriate if a simulation serves the sole purpose of finding an accurate approximation to the velocity field, but not necessarily of reflecting the exact sparsity pattern of the strain rate.

\subsection{Lid-Driven Cavity}

We now move on to the popular benchmark problem of viscoplastic flow inside a lid-driven square cavity. The fluid motion is only driven by a moving wall, i.e. we have $\bs f = 0$ and we define
\begin{equation*}
\bs u_\mathrm{D} (x) =
\begin{cases*}
(1,0)^\top & if $x_2 = 1$\\
(0,0)^\top & otherwise
\end{cases*}.
\end{equation*}
We point out that due to the discontinuities in the top corners, this choice violates the assumption $\bs u_\mathrm{D} \in U_\mathrm{D}$. The lid-driven cavity problem is classically studied to assess the performance of numerical methods under the presence of singularities.

\paragraph{Iterations and Computing Times}

For a range of different grid sizes and values of the Bingham number, we now compare how many iterations and how much time each of the four algorithms ADMM / ALG2, ISTA*, FISTA* and FISTA* with re-starting requires to compute a solution of the prescribed accuracy $\mathrm{gradTol} = 10^{-4}$. For these simulations, we focus on the problem of Bingham flow.

\begin{figure}
\centering
\includegraphics[width=\textwidth]{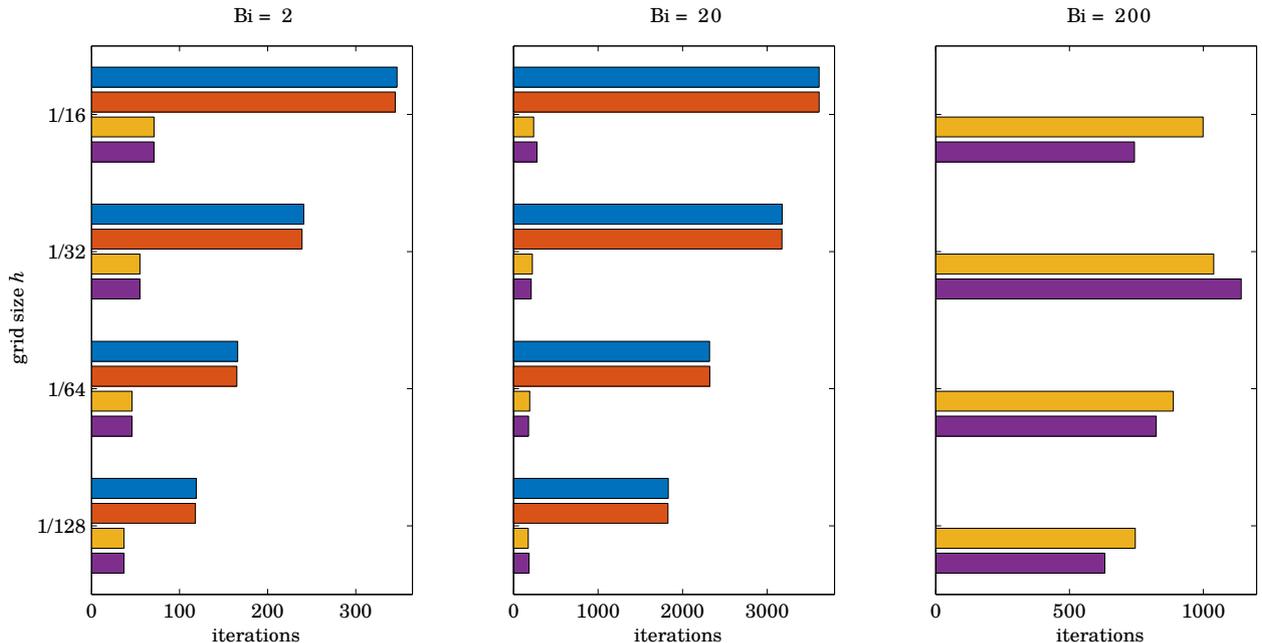}
\caption{Iterations until convergence for different mesh sizes (of the coarse pressure grid) and yield stress parameters. ALG2 (blue) and ISTA* (orange) failed to converge within 5,000 iterations for $\bi = 200$, therefore only the results for FISTA* (yellow) and FISTA* with re-starting (purple) are shown here.}
\label{fig:iter}
\end{figure}

\begin{figure}
\centering
\includegraphics[width=\textwidth]{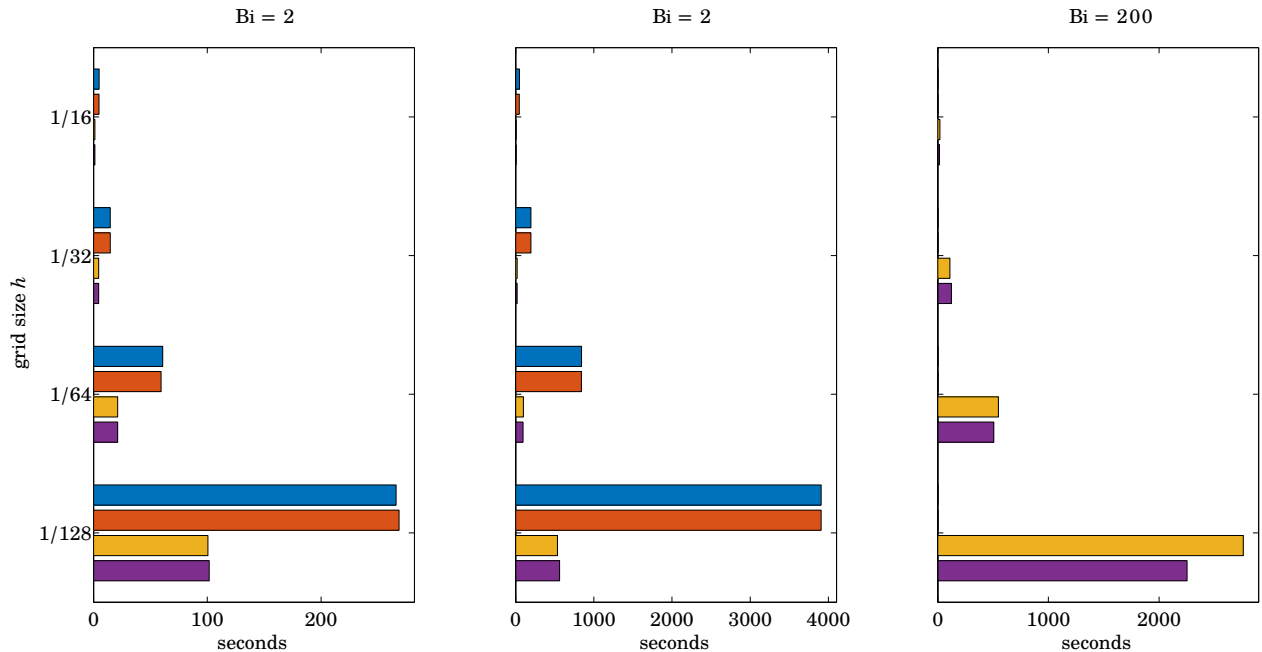}
\caption{Computing times corresponding to the test runs in Figure \ref{fig:iter}.}
\label{fig:time}
\end{figure}

It turns out again that by incorporating an acceleration scheme, the number of iteration is reduced significantly in every single case. For the largest value of the Bingham number considered here, $\bi = 200$, ALG2 and ISTA* were still far from an optimal solution even after 5,000 iterations, which we consider as \emph{failed to converge}. It is well possible that for different algorithm parameters ALG2 would have terminated within these 5,000 iterations. However, such parameters would have to be found by trial and error, which is clearly unpractical. For the sake of comparability with ISTA* and FISTA*, where such heuristics are not necessary, we limit our presentation to the setting $\varrho = s^{(k)} = 1/L^{(k)} \equiv 2$.

For the other cases, where we have data for all algorithms available, the dual FISTA method requires 83\% fewer iterations and 79\% less computing time than the alternating direction method of multipliers. The reduction in iteration numbers and CPU times for the re-starting adaptation are 83\% and 78\%, respectively. As can be seen from Figures \ref{fig:iter} and \ref{fig:time}, re-starting is worthwhile in certain examples, while it is the opposite in others.

\subsection{Yielded and Unyielded Flow Regions}

As noted by Yu and Wachs \cite{Yu2007}, the precise resolution of yielded and unyielded regions becomes particularly challenging in case of larger values of the yield stress.

Solutions for the problem, using different computational techniques and different values of the yield stress, have been published by Begis \cite{Begis1979}, (see also \cite[Ch~6]{Glowinski1989}), Sanchez \cite{Sanchez1998}, Mitsoulis and Zisis \cite{Mitsoulis2001}, Vola et al. \cite{Vola2003}, Yu and Wachs \cite{Yu2007}, Olshanskii \cite{Olshanskii2009}, De los Reyes and Gonz{\'a}lez Andrade \cite{Reyes2010}, Zhang \cite{Zhang2010}, Glowinski and Wachs \cite{Glowinski2011}, dos Santos et al. \cite{Santos2011}, Syrakos et al. \cite{Syrakos2013}, Aposporidis et al. \cite{Aposporidis2014} and Muravleva \cite{Muravleva2015}. For specific features of the flow, such as vortex positions and intensities, we also refer to these works.

\begin{figure}
\centering
\includegraphics[width=0.49\textwidth]{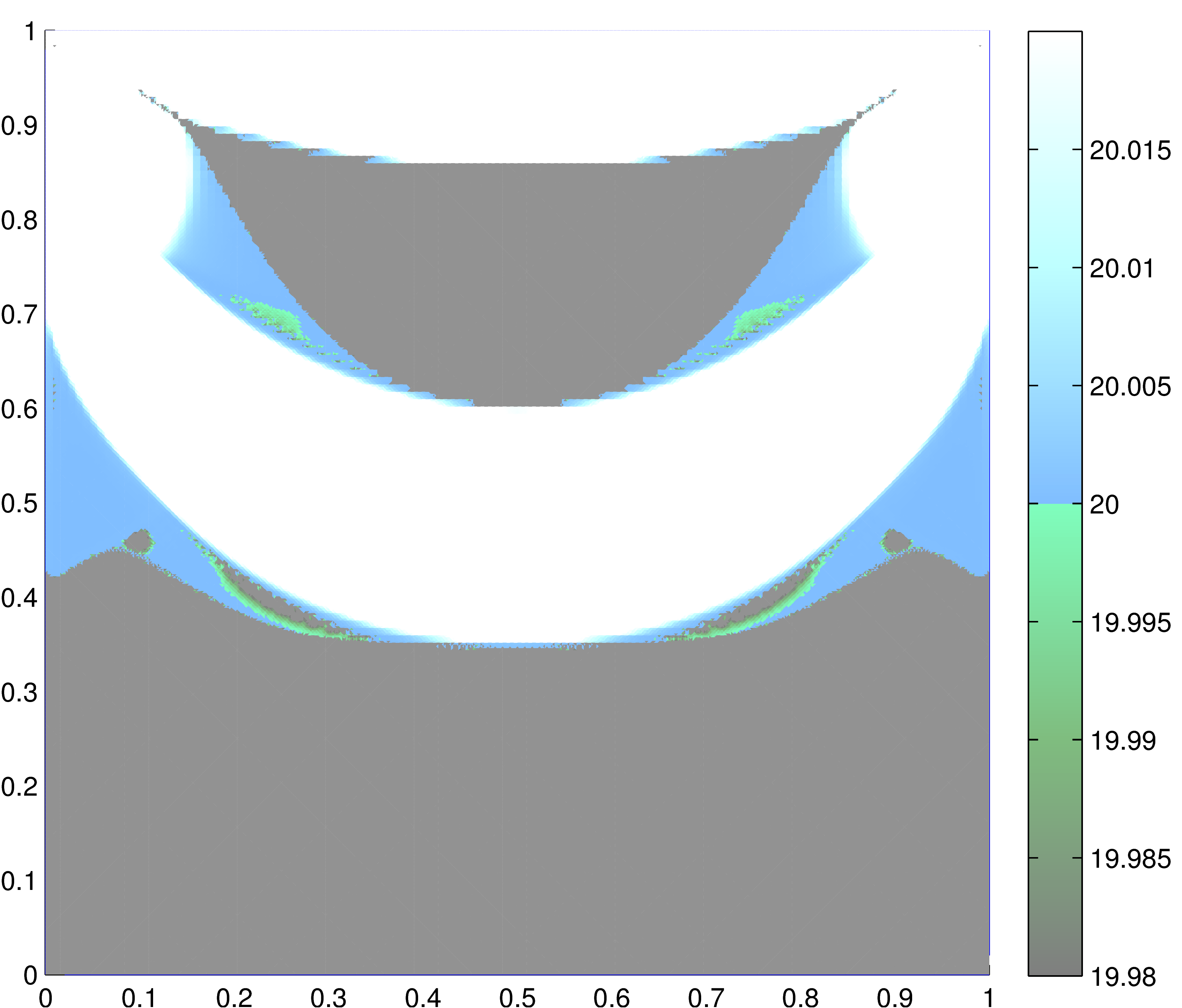}
\hfill
\includegraphics[width=0.49\textwidth]{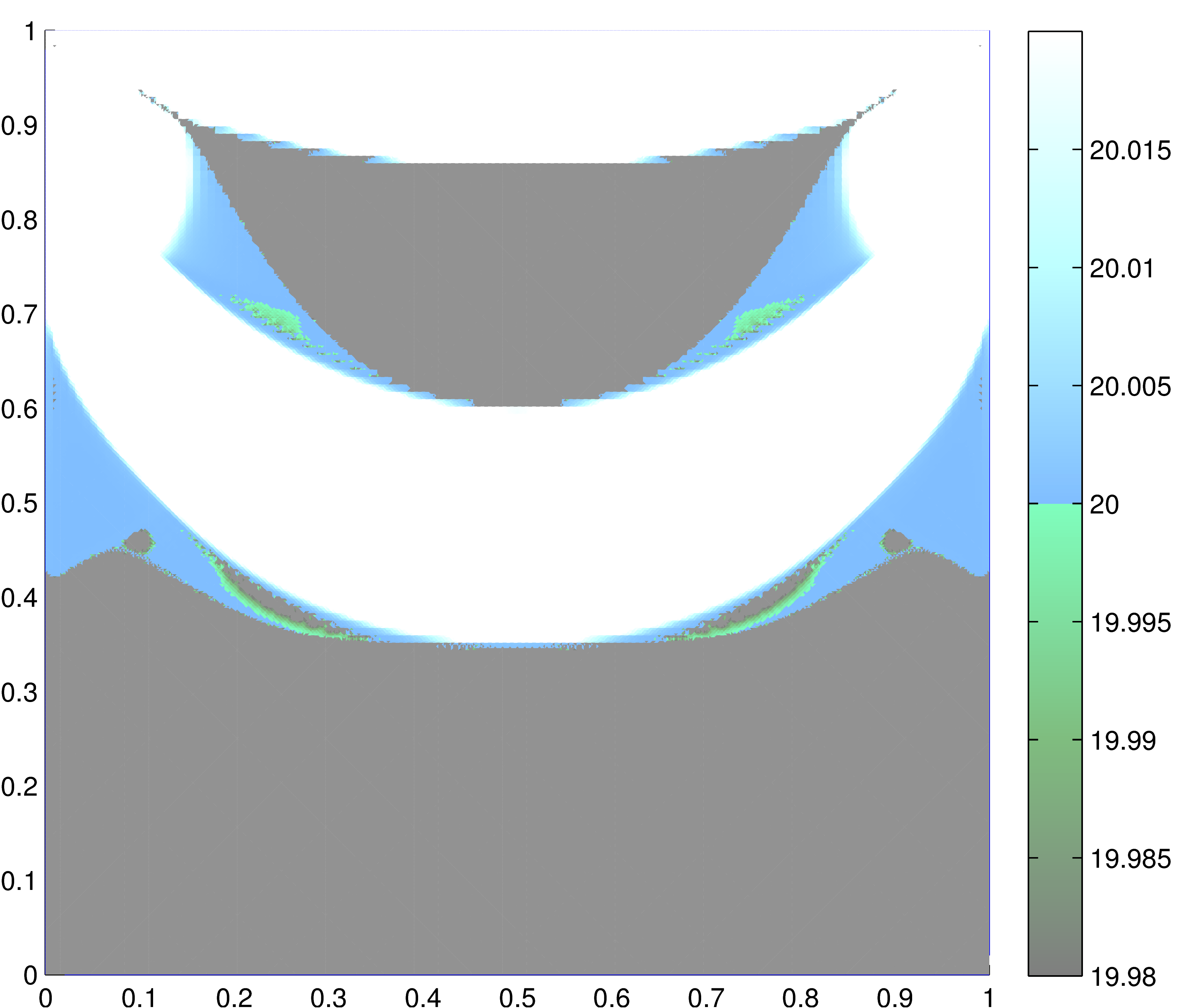}
\includegraphics[width=0.49\textwidth]{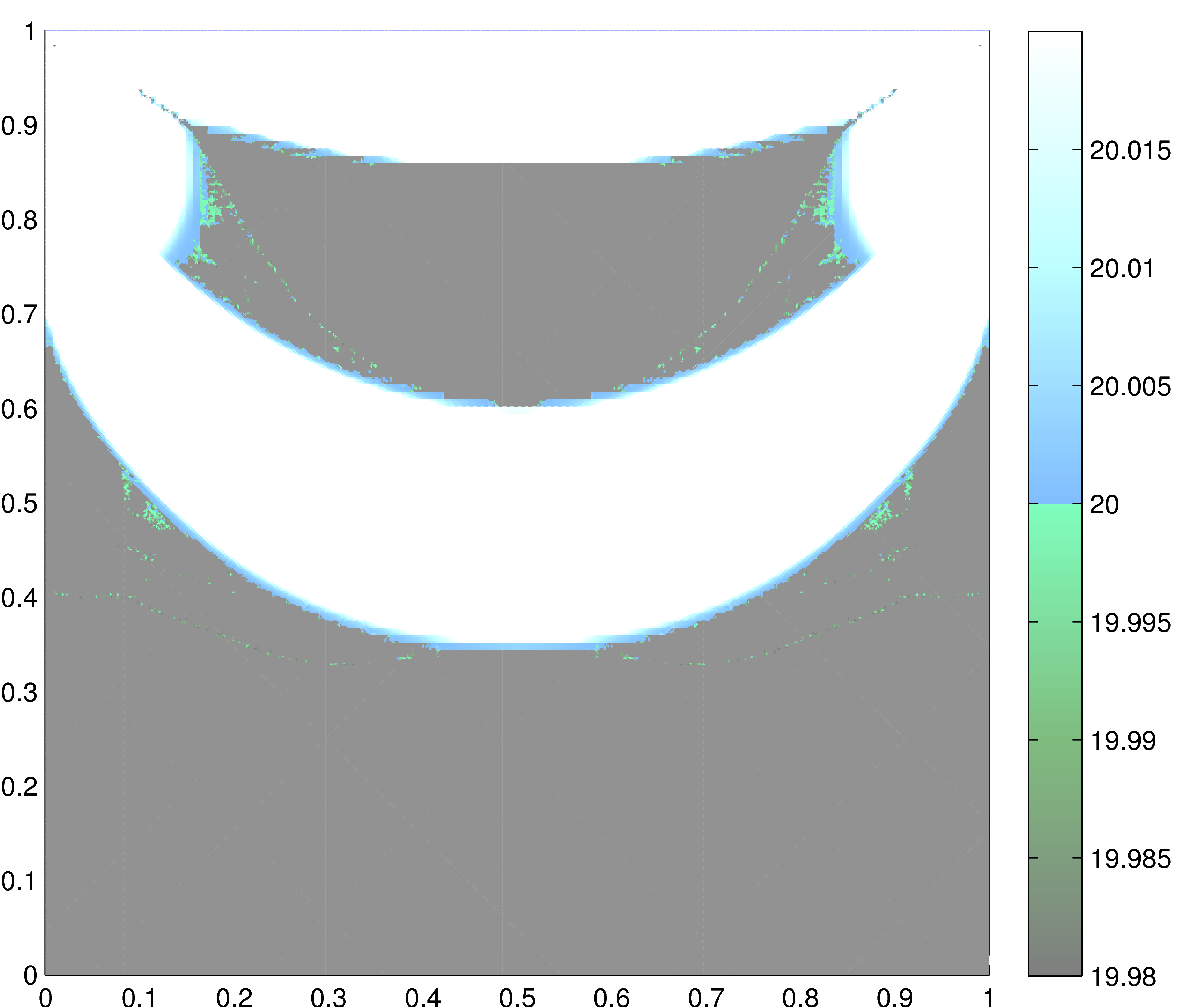}
\hfill
\includegraphics[width=0.49\textwidth]{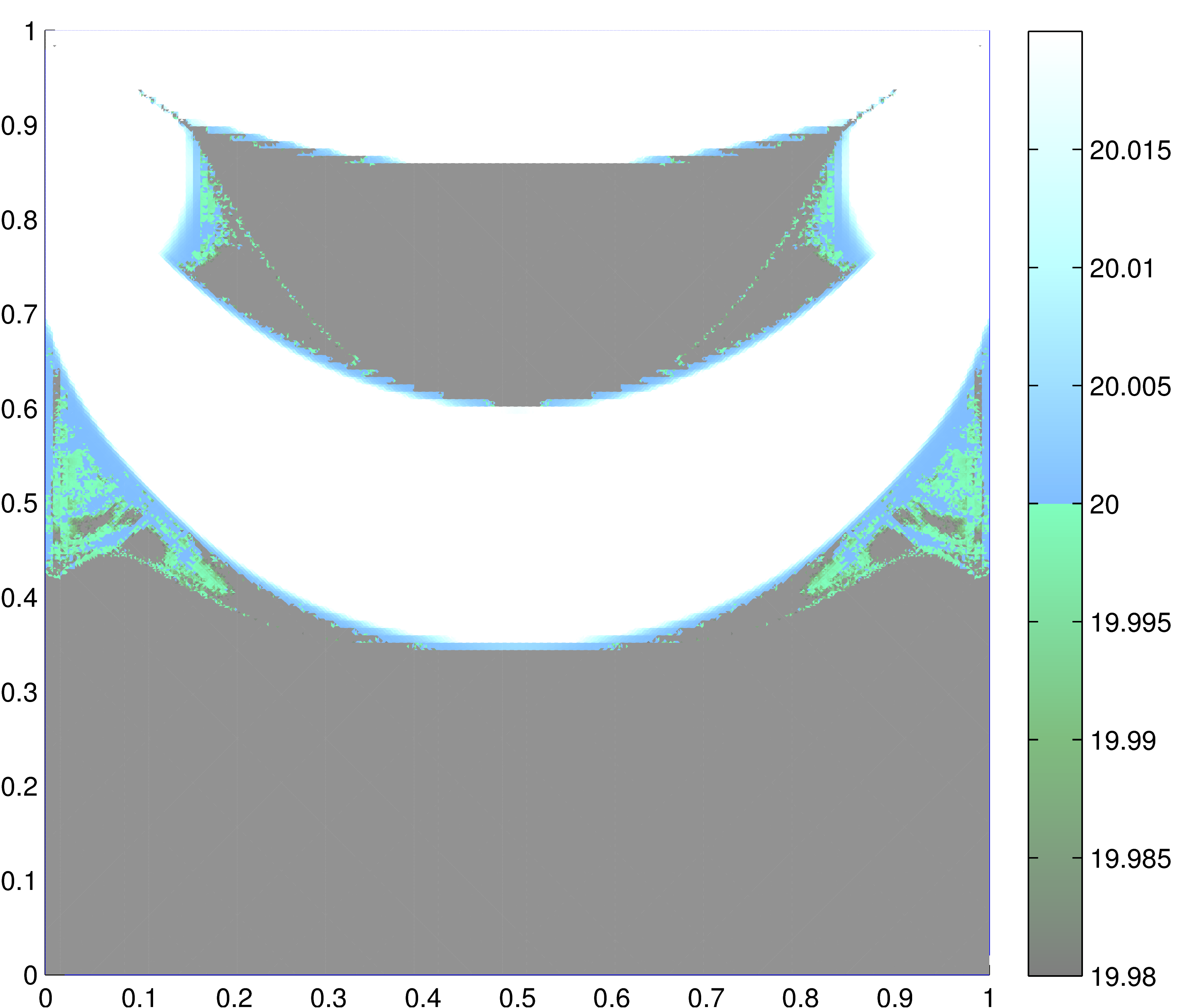}
\caption{\textsc{Frobenius} norm of the stress $\lvert \bs\tau \rvert$ near the Bingham number $\bi$. Top left to bottom right: ADMM, ISTA*, FISTA* and FISTA* with re-starts. Values outside the range of the colourbar have been projected onto the upper and lower end points, respectively ($\bi = 20$, $h = 1/128$).}
\label{fig:yield-surface-uniform}
\end{figure}

From Figure \ref{fig:yield-surface-uniform} we observe that Algorithm FISTA* identifies the unyielded regions in agreement with the results published in the works cited above. The approximation computed with FISTA* including re-starts is overall similar. Nevertheless, the relatively large areas where the stress is very close to the yield stress make it difficult to detect where the stagnant flow region ends and where shearing begins. Overall, elements where $\lvert \bs\tau_h \rvert \leq \bi$ clearly dominate in these areas, which should, indeed, be classified as unyielded.

Despite the identical stopping criterion in all cases, Algorithms ADMM and ISTA* clearly underestimate the regions occupied by unyielded fluid. While the approximation of the stress lies at least reasonably close to the yield stress in the blue areas, these two methods still fail to identify these as solid.

\paragraph{Model Reduction with Adaptive Finite Elements}

In past years, solutions on adaptive grids have already been successful at resolving the liquid-solid interface in fine detail, while reducing the substantial computational cost of simulations on uniform grids with the same fine resolution \cite{Zhang2010,Saramito2001,Roquet2003,Syrakos2014}. Similarly, our objective is to achieve a resolution of $h = 1/128$ in critical areas, while using a much coarser mesh with $h = 1/16$ where the residual is already comparatively small.

\begin{figure}
\centering
\includegraphics[width=0.49\textwidth]{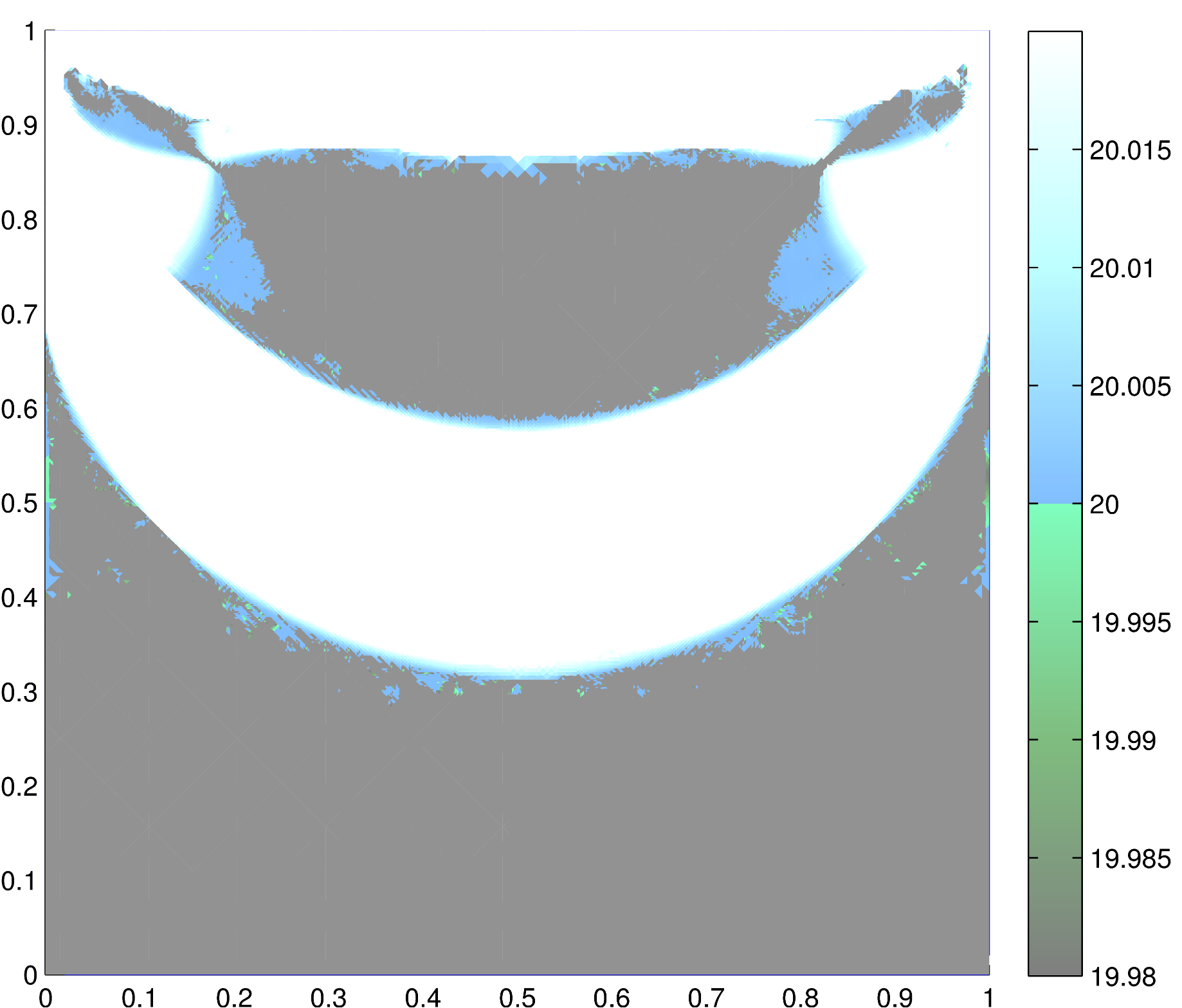}
\includegraphics[width=0.49\textwidth]{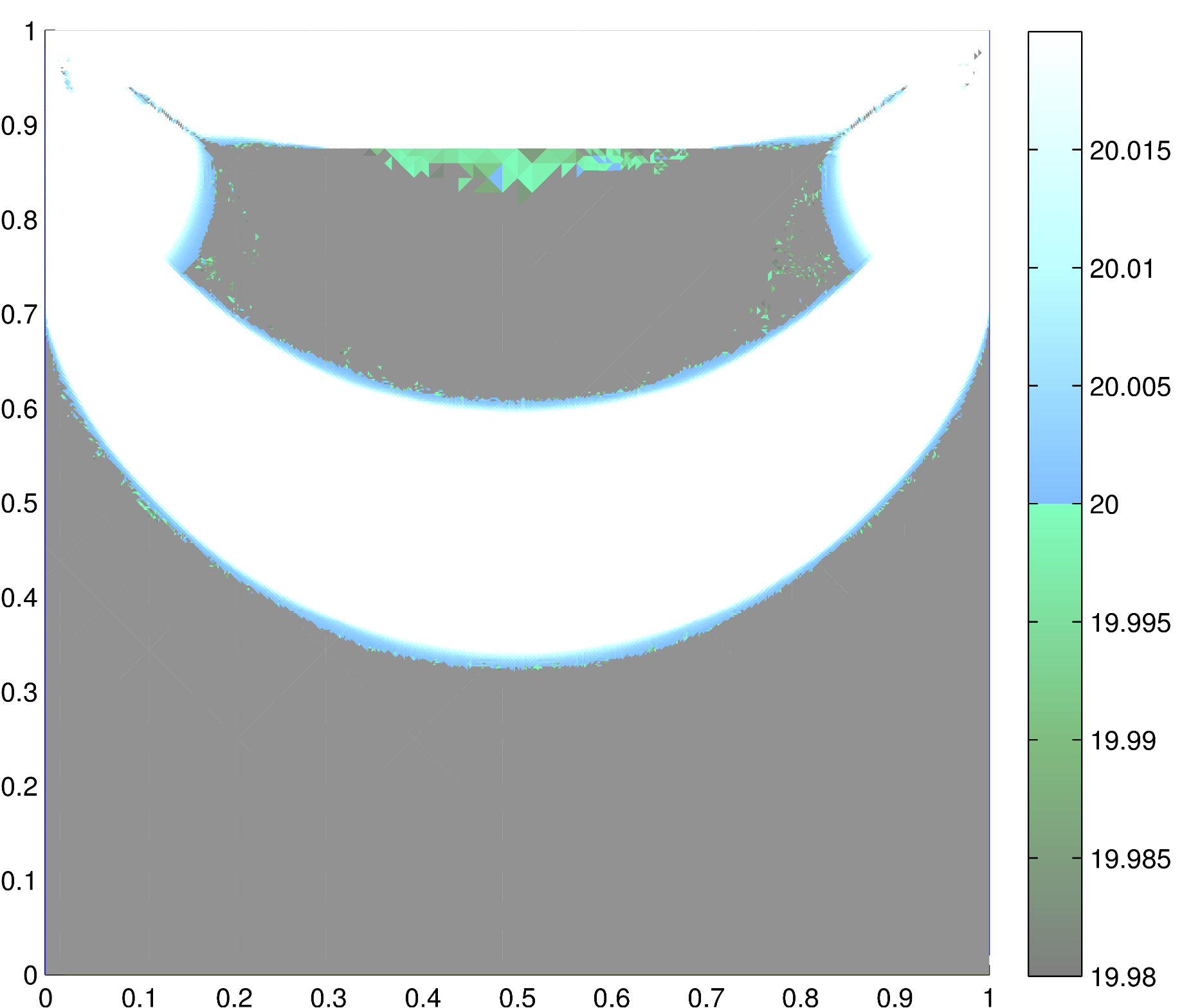}
\caption{Adaptive finite elements for resolving the free boundary between yielded and unyielded regions. Stress magnitude $\vert \bs \tau_h \vert$ computed by ALG2 (left) and FISTA* (right) for the same problem as in Figure \ref{fig:yield-surface-uniform} ($\bi = 20$).}
\label{fig:yield-surface-adaptive}
\end{figure}

For now, let us use the following ad hoc strategy:
\begin{itemize}
\item Solve the optimisation problem with one of the four algorithms until convergence.
\item Determine the 60th percentile of the Frobenius norm of the residual $\lvert\sg\bs u^{(k)}_h - \bs \sr^{(k)}_h\rvert$ over all triangles and refine those $\sim 40\%$ of all triangles with the largest residual. Further refinements of neighbouring triangles are required to avoid hanging nodes.
\item Interpolate the converged solution linearly to the refined grid.
\end{itemize}

In Figure \ref{fig:yield-surface-adaptive}, we tackled the problem with $\bi = 20$ with ALG2 and FISTA* once again, this time on a grid that was only locally refined. Starting from the uniform mesh with $h = 1/16$, we cycled through the above refinement procedure three times. We conclude that the quality of both results is very much comparable to the one of the corresponding graphs in Figure \ref{fig:yield-surface-uniform}. Nevertheless, it took about 65\% (ALG2) or 61\% (FISTA*) less computing time, respectively, until convergence was achieved. Additionally, the identification of the zero-flow region by ALG2 has even improved considerably. The upper stagnant zone still exhibits many coarse artefacts, though.

\begin{figure}
\centering
\includegraphics[width=0.53\textwidth]{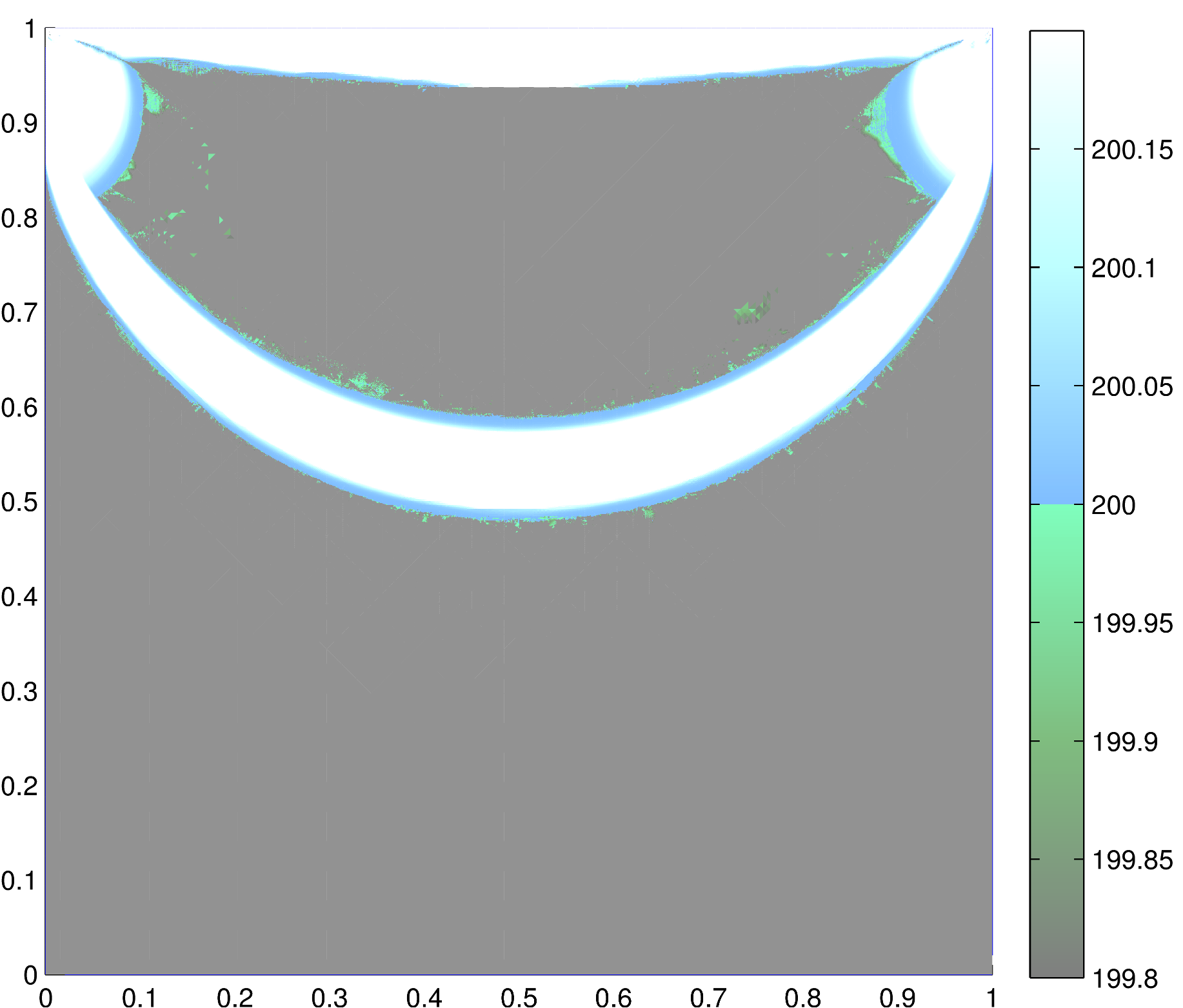}
\includegraphics[width=0.45\textwidth]{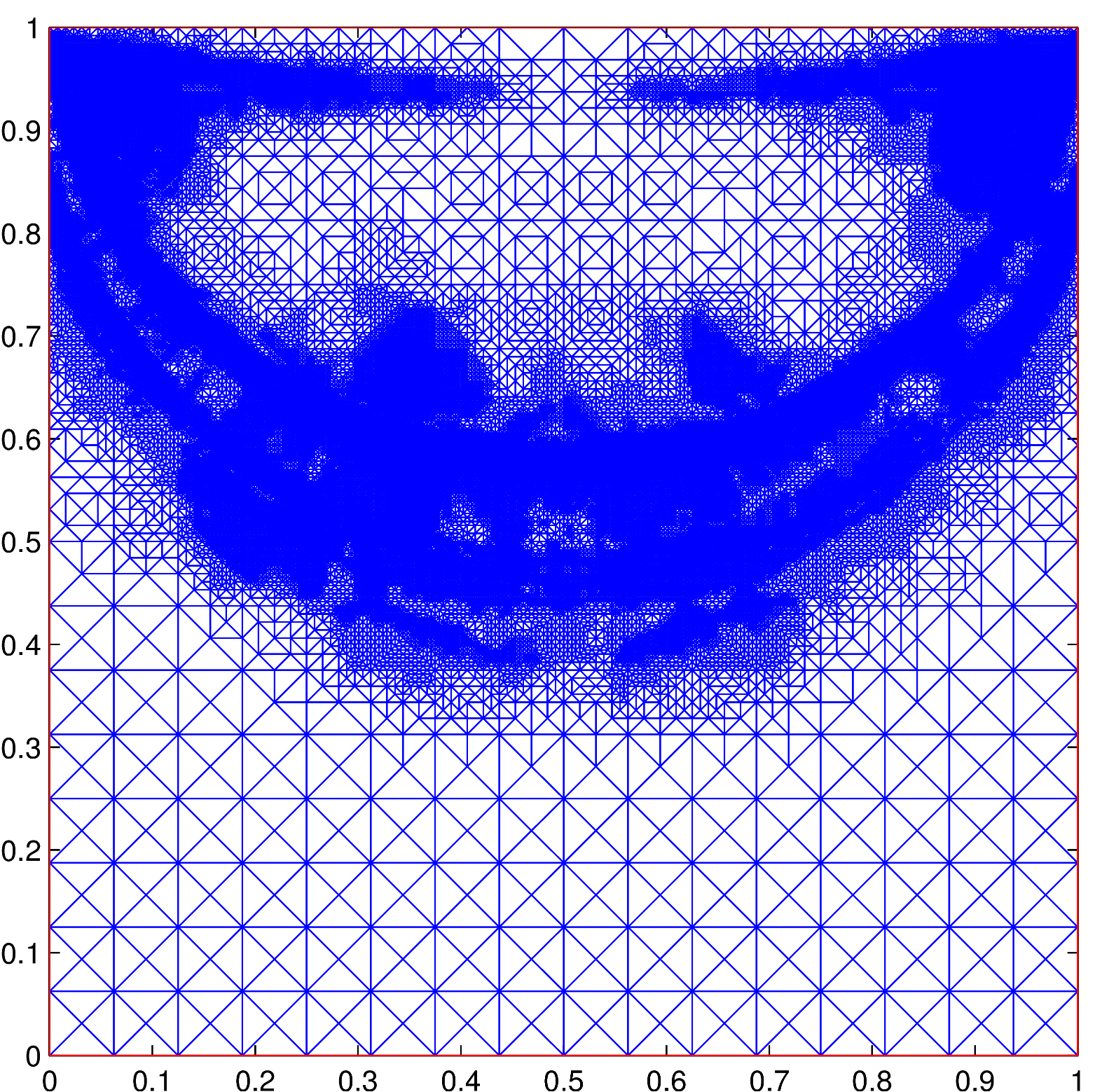}
\includegraphics[width=0.53\textwidth]{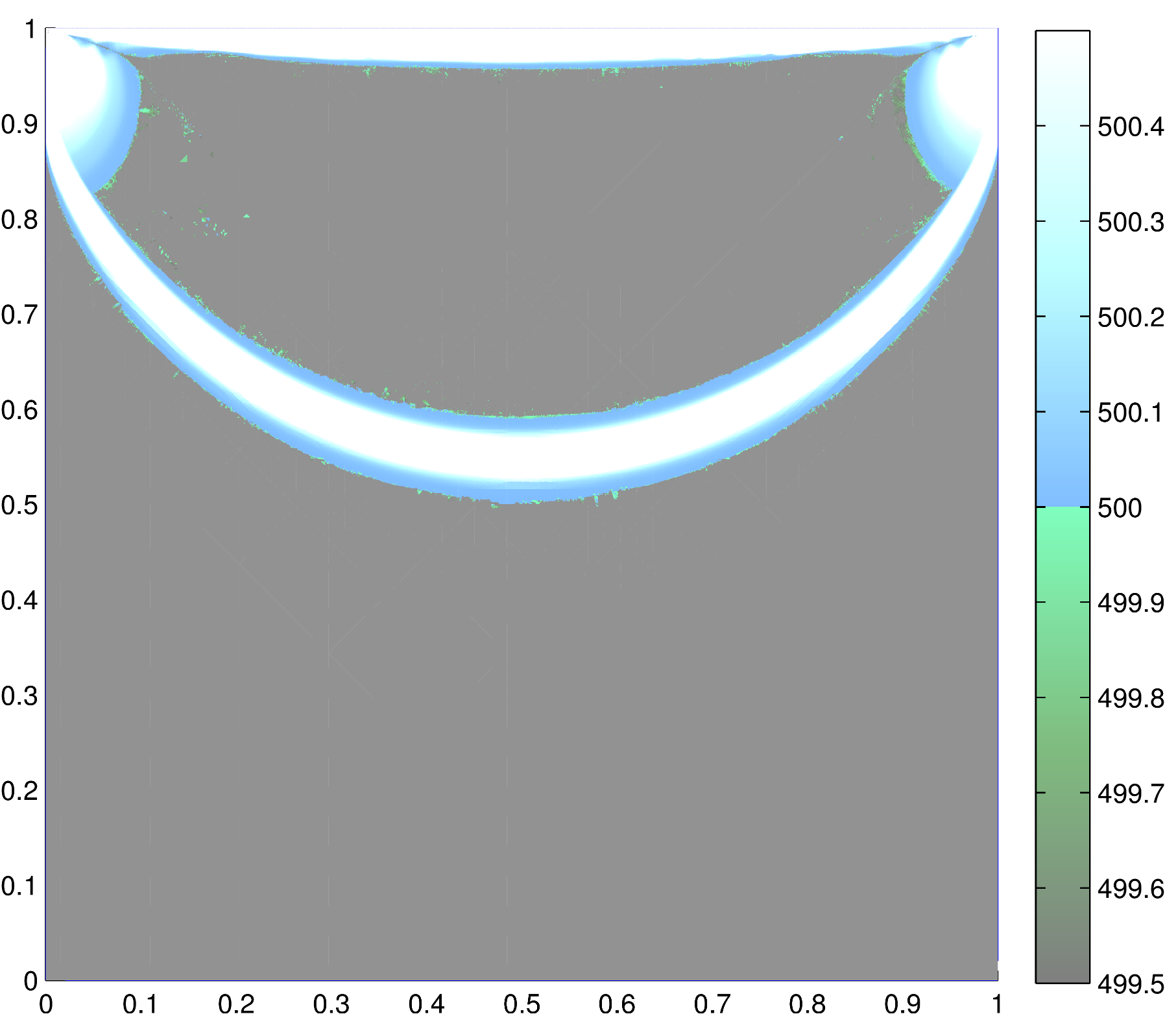}
\includegraphics[width=0.45\textwidth]{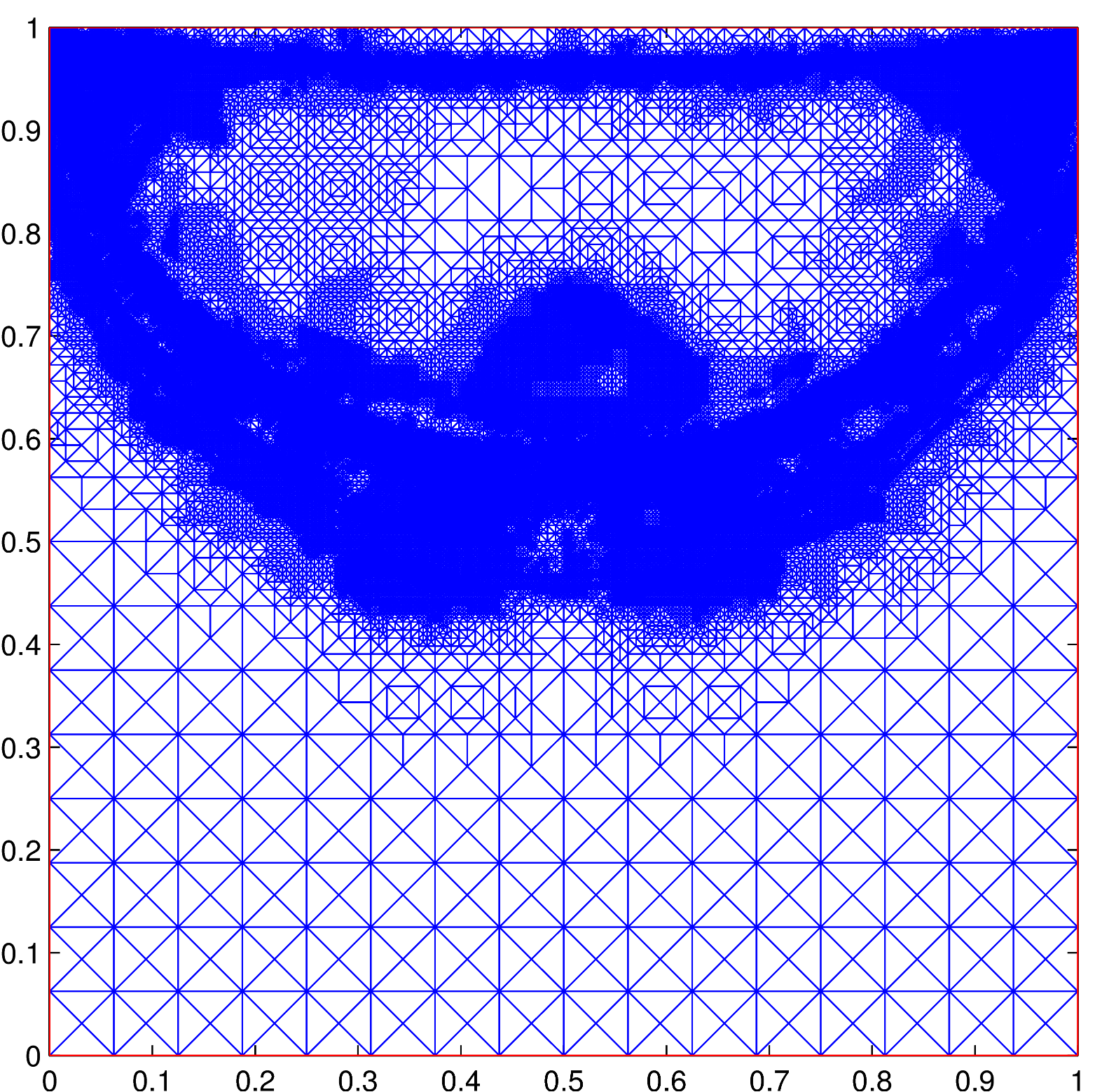}
\caption{Predicting the yield surface for flows at high yield stress values with FISTA* and adaptive finite elements. Top: $\vert \bs \tau_h \vert$ and the corresponding grid when $\bi = 200$. Bottom: $\bi = 500$.}
\label{fig:yield-surface-200-500}
\end{figure}

Now that we can assume our refinement methodology to be validated for FISTA*, we apply it to the traditional challenge of predicting the yield surfaces when the yield stress is very large. We pick the two values $\bi = 200$ and $\bi = 500$ for which results have been published in the literature. We adaptively refine the initial homogeneous mesh with $h = 1/16$ five times. Our results are depicted in Figure \ref{fig:yield-surface-200-500}.

Yu and Wachs \cite{Yu2007} and Muravleva \cite{Muravleva2015} have used ALG2 on a homogeneous grid with $h = 1/256$ to solve these two problems. Their results deviate from each other as well as from ours, which are in close qualitative agreement with the publications of Mitsoulis and Zisis \cite{Mitsoulis2001} and Syrakos et al. \cite{Syrakos2013}. Since both of the latter works solve regularised approximations of the Bingham flow problem, fine geometric features like sharp tips that are visible in our results, have already been smoothed out in the problem formulations of these authors.

Although this very basic approach to mesh adaptivity has proven to be effective, we anticipate even further improvements from more sophisticated, goal-oriented adaptive finite element methods like the DWR (dual weighted residual) method. We refer to the book of Suttmeier \cite{Suttmeier2008} for more details.

\subsection{Conclusions}

We wish to emphasise that unlike second-order methods for instance, which require a Hessian at every iteration, the higher rate of convergence of the accelerated dual proximal gradient method FISTA* compared to the classical alternating direction method of multipliers ADMM / ALG2 comes at the very minimal cost of additionally evaluating a linear combination and storing one extra variable. Even though some open questions remain regarding the convergence of both FISTA* and ALG2, these seem to be of no practical relevance. For FISTA*, we have presented strategies related to the definition of the primal sequences, which allow us to generally prove an accelerated worst-case convergence rate of $O(1/k)$, compared to only $O(1/\sqrt{k})$ for classical algorithms.

Furthermore, globally optimal values for free parameters that occur in the dual proximal gradient methods can either be calculated a priori, or estimated numerically by backtracking. Moreover, the new dual FISTA method is very closely related to ALG2, in the sense that the subproblems that arise in both algorithms are either identical or even simpler for FISTA*. Any existing code based on an augmented Lagrangian method can therefore easily be modified to implement FISTA*. This is what leads us to our conclusion that the Algorithm FISTA* could be seen as a more efficient successor algorithm of ALG2 for solving genuinely nonsmooth formulations of viscoplastic flow problems accurately.

We believe that at this stage it is still too early to express a recommendation towards either the simple FISTA* method or the variant with adaptive re-starting. Our small number of examples so far indicate that re-starting may not be be quite as effective in the context of Bingham flow as it is for other nonsmooth optimisation problems \cite{ODonoghue2013}. Further numerical studies are required to provide more guidance on this question.

Our analysis already applies to a more general framework, including viscoplastic flow in three spatial dimensions. The extension to time-dependent flow problems follows naturally by first applying a suitable semi-discretisation scheme in time. If the inertial term is is discretised explicitly, then, at every time step, the methodology of this paper remains applicable except that an additional mass matrix arises in the solution of the Stokes problems. This strategy is completely analogous to solutions by other numerical methods, e.g. \cite{Reyes2012,Muravleva2015}.

We still see potential for improving the efficiency of FISTA* further, e.g. by employing preconditioning techniques (cf \cite{Bonettini2015}) or inexact evaluations of the proximal map (cf \cite{Schmidt2011,Jiang2012}). These concepts shall be our focus of further research on the topic.

\section*{Acknowledgements}

The authors wish to thank Shoham Sabach from Technion in Haifa, Israel, for stimulating this fruitful investigation of accelerated gradient schemes in viscoplasticity.

\bibliographystyle{elsarticle-num} 
\bibliography{bibliography}

\end{document}